\documentclass[letterpaper]{amsart}

\usepackage{amssymb}
\usepackage{amsthm}
\usepackage{amsmath}
\usepackage{bm}
\usepackage{amsfonts}
\usepackage{yfonts}
\usepackage{comment}
\usepackage[mathscr]{euscript}
\usepackage{upgreek}
\usepackage[T1]{fontenc}
\usepackage{ wasysym }

\usepackage{stackengine}
\stackMath
\usepackage{tikz-cd}
\usepackage[all]{xy}
\usepackage{xfrac}
\usepackage{MnSymbol}
    
\usepackage{soul}

\usepackage{todonotes}
\usepackage{xargs}

\newcommand{\newdownfree}{\scalebox{1.5}[1]{\ensuremath{\downfree}}}
\newcommand{\newndownfree}{\scalebox{1.5}[1]{\ensuremath{\ndownfree}}}

\newcommandx{\margin}[2][1=]{\todo[linecolor=blue,backgroundcolor=blue!25,bordercolor=blue,#1]{#2}}

\def\Ind#1#2{#1\setbox0=\hbox{$#1x$} \kern\wd0    \hbox to 0pt{\hss$#1\newdownfree$\hss}\kern\wd0}
\def\Notind#1#2{#1\setbox0=\hbox{$#1x$}\kern\wd0 \hbox to 0pt{\hss$#1\newndownfree$\hss}\kern\wd0}

\newcommand{\lar}{L_\mathrm{f}}
\newcommand{\lart}{L^\times_\mathrm{f}}
\newcommand{\lam}{L_\mathrm{mc}}
\newcommand{\lag}{L_\mathrm{m}}
\newcommand{\laa}{L_\mathrm{al}}
\newcommand{\lac}{L_\mathrm{t}}

\newcommand{\thmp}{\mathrm{GMO}_p}
\newcommand{\thmm}{\mathrm{GMO}}
\newcommand{\thmo}{\mathrm{GMO}_0}
\newcommand{\thap}{\mathrm{HAO}_p}
\newcommand{\thaa}{\mathrm{HAO}}
\newcommand{\thao}{\mathrm{HAO}_0}
\newcommand{\thgg}{\mathrm{GM}}
\newcommand{\thgp}{\mathrm{GM}_p}
\newcommand{\thgo}{\mathrm{GM}_0}

\newcommand{\acl}{\mathrm{acl}}

\newcommand{\elesub}{\preccurlyeq}

\newcommand{\qq}{\mathbb{Q}}
\newcommand{\bbT}{\mathbb{T}}
\newcommand{\rr}{\mathbb{R}}
\newcommand{\ff}{\mathbb{F}}
\newcommand{\nn}{\mathbb{N}}
\newcommand{\uu}{\mathbb{U}}
\newcommand{\zz}{\mathbb{Z}}

\newcommand{\cc}{\mathbb{C}}

\newcommand{\la}{\langle}
\newcommand{\ra}{\rangle}

\newcommand{\ACF}{\mathrm{ACF}}
\newcommand{\NTP}{\mathrm{NTP}}
\newcommand{\ACFO}{\mathrm{ACFO}}

\newcommand{\sM}{\mathscr{M}}

\newcommand{\sP}{\mathscr{P}}

\DeclareMathOperator{\mult}{mult}

\makeatletter
\newsavebox\myboxA
\newsavebox\myboxB
\newlength\mylenA

\newcommand*\xbar[2][0.75]{%
    \sbox{\myboxA}{$\m@th#2$}%
    \setbox\myboxB\null
    \ht\myboxB=\ht\myboxA%
    \dp\myboxB=\dp\myboxA%
    \wd\myboxB=#1\wd\myboxA
    \sbox\myboxB{$\m@th\overline{\copy\myboxB}$}
    \setlength\mylenA{\the\wd\myboxA}
    \addtolength\mylenA{-\the\wd\myboxB}%
    \ifdim\wd\myboxB<\wd\myboxA%
       \rlap{\hskip 0.5\mylenA\usebox\myboxB}{\usebox\myboxA}%
    \else
        \hskip -0.5\mylenA\rlap{\usebox\myboxA}{\hskip 0.5\mylenA\usebox\myboxB}%
    \fi}
\makeatother

\newcommand{\ir}{(i)}
\newcommand{\iir}{(ii)}
\newcommand{\iiir}{(iii)}
\newcommand{\ivr}{(iv)}

\newcommand{\meno}{\medskip \noindent}

\newtheorem{thm}{Theorem}[section]
\newtheorem{lem}[thm]{Lemma}
\newtheorem*{thm*}{Theorem}
\newtheorem*{thmA*}{Theorem A}
\newtheorem*{prop*}{Proposition}
\newtheorem*{cor*}{Corollary}
\newtheorem*{corB*}{Corollary B}
\newtheorem*{corC*}{Corollary C}
\newtheorem{prop}[thm]{Proposition}
\newtheorem{cor}[thm]{Corollary}
\theoremstyle{definition}

\newtheorem{rem}[thm]{Remark}
\newtheorem{fact}[thm]{Fact}

\definecolor{red}{rgb}{1.0, 0, 0}

\begin{document}
\sloppy

\title[Tame structures via character sums over finite fields]{Tame structures via character sums \\ over finite fields}
\author{Chieu-Minh Tran}
\address{Department of Mathematics, University of Illinois at Urbana-
Champaign, Urbana, IL 61801, U.S.A}
\curraddr{}
\email{mctran2@illinois.edu}
\date{\today}

\begin{abstract}
We show that the theory of  algebraically closed fields with multiplicative circular orders has a model companion $\ACFO$. Using number-theoretic results on character sums over finite fields, we show that if $\ff$ is an algebraic closure of a finite field, and $\triangleleft$ is any translation-invariant  circular order on the multiplicative group $\ff^\times$, then $(\ff, \triangleleft)$ is a model of $\ACFO$. Our results can be regarded as analogues of Ax's results in \cite{Ax} which utilize counting points over finite fields.

\end{abstract}

\maketitle

\section{Introduction}

\noindent Throughout,  $\ff$ will be an algebraic closure of a finite field. We are interested in the following question:
\begin{center}
    {\it Are there natural expansions of $\ff$ by order-type relations which are also model-theoretically tame?}
\end{center}
There is no known order-type relation on $\ff$ which interacts in a sensible way with both addition and multiplication. This is in stark contrast to the situation with the field $\cc$ of complex numbers where addition and multiplication are compatible with the  Euclidean metric induced by the natural order on $\rr$. It is not hard to see the reason: the additive group of $\ff$ is an infinite torsion group of finite exponent, so even finding an additively compatible  order-type relation seems unlikely. On the other hand, the multiplicative group $\ff^\times$ is a union of cyclic groups, so it is fairly natural to consider  circular orders $\triangleleft$ on $\ff^\times$ which are compatible with the multiplicative structure. In this paper, we will show that the resulting structures $(\ff, \triangleleft)$ give a positive answer to some aspects of the above question.

\meno We will take a step back to be more precise and to study the above structures as members of a natural class.  A {\bf circular order} on a group $G$ is a ternary relation $\triangleleft$ on $G$ which is invariant under multiplication by elements in $G$ and satisfies the following conditions for all $a,b,c \in G$:
\begin{enumerate}
\item if $\triangleleft(a,b,c)$, then $\triangleleft (b,c,a) $;
\item if $\triangleleft (a,b,c)$, then not $\triangleleft (c,b,a)$;
\item if $\triangleleft (a,b,c)$ and $\triangleleft(a,c,d)$, then $\triangleleft(a,b,d)$;
\item if $a,b,c$ are distinct, then either $\triangleleft(a,b,c) $ or $\triangleleft(c,b,a) $.
\end{enumerate}
A canonical example, also used later on, is $(\bbT, \triangleleft)$ where $\bbT$
is the multiplicative group of complex numbers with norm $1$, and $\triangleleft$ is the  clockwise circular order (i.e., $\triangleleft(a,b,c)$ if $b$ lies in the clockwise open arc from  $a$ to $c$ viewing $\bbT$ as the unit circle).

\meno A \textbf{multiplicative circular order} on a field $F$ is a circular order on the multiplicative group $F^\times$, viewed as a ternary relation on $F$. If $\triangleleft$ is a multiplicative circular order on a field $F$, then $(F, \triangleleft$) is a structure in the total language $\lac$ extending the language $\lar = \{ 0, 1, +, -, \times, \square^{-1}\}$ of fields by a ternary predicate symbol $\triangleleft$. Let $\ACFO^-$ be the $\lac$-theory whose models are such $(F, \triangleleft)$ where $F$ is algebraically closed. Section~\ref{sec: Model  companion of GMO-} and Section~ \ref{sec: Model  companion of ACFO-} establish our first main result:
 
\begin{thm} \label{thm: Model companion of ACFO-}
The theory $\ACFO^-$ has a model companion $\ACFO$.
\end{thm}
\noindent  Underlying the proof of Theorem~\ref{thm: Model companion of ACFO-} is the following heuristic: the existential closed models of $\ACFO^-$ are $(F, \triangleleft) \models \ACFO^-$ where  $(F^\times, \triangleleft)$ is ``sufficiently rich'', and  $+$ interacts in a ``random fashion'' with $\triangleleft$ modulo their compatibility with $\times$. The challenges involve making sense of  ``sufficiently rich'' and ``random fashion'',  justifying this heuristic, and showing that these properties are first-order axiomatizable. 

\meno We  return to the structures described in the first paragraph in Section~\ref{Sec: Standard models are existentially closed}: 
\begin{thm} \label{thm: genericity of the standard models}
If $\triangleleft$ is a multiplicative circular order on $\ff$, then
$(\ff, \triangleleft)\models \ACFO$.
\end{thm}
\noindent   Every injective group homomorphism  $\chi: \ff^\times \to \mathbb{T} $ induces a  multiplicative circular order on $\ff$, namely, the pullback $\triangleleft_\chi$ of  the clockwise circular order $\triangleleft$ on $\mathbb{T}$ by the map $\chi$. It turns out that every multiplicative circular order on $\ff$ is of this form; see Corollary~\ref{Up is model of GMO-3}. The main idea of the proof of Theorem~\ref{thm: genericity of the standard models} is to exploit this connection and results on character sums over finite fields. These results are useful here as they reflect ``number-theoretic randomness''~\cite{Konyagin}. This is precisely what we want for the interaction between $+$ and $\triangleleft$.  
 
\meno This work is a response to the question below by van den Dries and Hrushovski;  Kowalski also asked a related question in \cite{Kowalski}.

 \begin{center}
     {\it Do results on exponential sums and character sums over finite fields yield any model-theoretically tame structures?}
 \end{center}
Behind this question is the hope to find  analogies of Ax's results in \cite{Ax}. There, the model-theoretic tameness of ultraproducts of finite fields essentially follows from results on counting points over finite fields. The theory $\ACFO$ is our proposed counterpart of the theory of pseudo-finite fields. Theorem~\ref{thm: Model companion of ACFO-} corresponds to the fact that the theory of pseudo-finite fields is almost model complete. Theorem~\ref{thm: genericity of the standard models} is an analogue of the fact that nonprincipal ultraproducts of finite fields are pseudo-finite (in the  definition given by Ax). There are also reasons to believe that there are deeper connections between $\ACFO$ and  the theory of pseudo-finite fields. Both theories include certain ``random features'' and can be put under the framework of interpolative fusions \cite{Interpolativefusion}, a subsequent joint work with Kruckman and Walsberg.

\meno Our earlier write-up~\cite{Version1} contained more results. Some of these are now generalized into results about interpolative fusions \cite{Interpolativefusion}. We do not include them here to minimize overlapping. The results on decidability are of a slightly different flavor, so they will be presented in a future paper.

\meno The structures $(\ff, \triangleleft)$ in Theorem~\ref{thm: genericity of the standard models} are not simple (in the sense of model theory) as they define dense linear orders. They also have $\text{IP}$ by a result of Shelah and Simon \cite{ShelahSimon}. It turns out that these structures do not even have $\text{TP}_2$ (see Proposition~\ref{Prop: this structure has TP_2}). This brings $(\ff, \triangleleft)$ outside the current known boundary of the combinatorially tame universe. We hope these structures will provide some motivation to push the boundary further and include them as well.
 
\subsection*{Notations and conventions} We keep the notational conventions of the introduction.
Moreover, we assume that $k$ and $l$ range over $\zz$, $m$ and $n$ range over $\nn=\{0, 1, 2,  \ldots\}$, $x= (x_1, \ldots, x_m)$ is an $m$-tuple of variables, $y=(y_1, \ldots, y_n)$ is an $n$-tuple of variables, $p$ is a prime, $q =p^k$ for some $k\geq 1$, $F$ is a field, $F^\times$ is the multiplicative group of $F$.

\section{Almost model companion of $\thmm^-$} \label{sec: Model  companion of GMO-} 

\noindent For understanding $\ACFO^-$ and finding its model companion, we need to first understand $F^\times$ and $(F^\times, \triangleleft)$ as $(F, \triangleleft)$ ranges over the models of $\ACFO^-$. Two phenomena turn out to be important later on:
\begin{enumerate}
    \item the ``reduct'' of $\ACFO^-$ to the language of multiplicative groups  is very simple: it ``almost'' admits quantifier elimination and has a natural notion of dimension;
    \item the ``reduct'' of  $\ACFO^-$ to the language of circularly ordered multiplicative groups ``almost'' has a model companion.
\end{enumerate}

\subsection{Multiplicative groups of algebraically closed fields} \label{Sec: Multiplicative groups of algebraically closed field} 
 
We will consider the theory of multiplicative groups of algebraically closed fields (i.e., the set of statements which hold in all such structures) in a suitable language, show that this theory ``almost'' admits quantifier elimination  and coincides with the theory of multiplicative groups of $\ACFO^-$-models, and obtain an axiomatization along the way as usual.

\meno Throughout Section~\ref{Sec: Multiplicative groups of algebraically closed field}, $G$ is a multiplicative abelian group. If $a$ and $b$ in $G$ are such that $b^n =a$, we call $b$ an {\bf $n$th root} of  $a$. If $n\geq 1$, an $n$th root of the identity element $1_G$ of  $G$ is {\bf trivial} if it is $1_G$. An $n$th root of $1_\mathbb{T}$ in $\mathbb{T}$ for some $n \geq 1$ is called a {\bf root of unity}. Let  $\uu \subseteq \mathbb{T}$ be the multiplicative group of roots of unity. For a given $p$, let  $\uu_{(p)} \subseteq \mathbb{T}$ be the multiplicative group of roots of unity whose order is coprime to $p$. So $\uu_{(p)}$ is isomorphic to $\ff^\times$ as a group when $\text{char}(\ff) =p$.

\meno Let $\lag = \{1, \times, \square^{-1} \}$ be the language of multiplicative groups.  It is easy to obtain an $\lag$-theory $\thgg$ such that $G \models \thgg$ if and only if the following conditions hold:
\begin{enumerate}
    \item[(G0)]  every finite subgroup of $G$ is cyclic;
    \item[(G1)]  the group $G$ is divisible;
    \item[(G2)]  for any two distinct prime numbers $p$ and $l$, either $1_G$ has a nontrivial $p$th root, or $1_G$ has a nontrivial $l$th root. 
\end{enumerate}
The theory  $\thgg$ is our candidate for axiomatizing the theory of multiplicative groups of algebraically closed fields. It has several natural extensions.
For a given $p$, let $\thgp$ be the $\lag$-theory whose models are $G \models \thgg$ which satisfy the following extra property:
\begin{enumerate}
    \item[(GC$_{p}$)]  every $p$th root of $1_G$ is trivial.
\end{enumerate}
It is easy to see that if $G \models \thgp$, then $1_G$ has a nontrivial $l$th root for any prime number $l \neq p$.
Let $\thgo$ be the $\lag$-theory whose models are  $G \models \thgg$ which satisfy the following extra property:
\begin{enumerate}
    \item[(GC$_{0}$)]  for all  prime numbers $l$,  $1_G$ has a nontrivial $l$th root.
\end{enumerate}
Hence,  a model of $\thgg$ is either a model of $\thgp$ for some $p$ or a model of $\thgo$.
\begin{rem} \label{Rem: the remark on roots}
Suppose $G$ satisfies conditions (G0) and (G1), and $l$ is a prime number. The condition that $1_G$ has a nontrivial $l$-root is also equivalent to several other conditions:
\begin{enumerate}
    \item $1_G$ has exactly $l$ many $l$th roots;
    \item for all $k\geq 1$, $1_G$ has exactly $l^k$ many $l^k$th  roots;
    \item for all $k\geq 1$, every $a \in G$ has exactly $l^k$th  many $l^k$th  roots.
\end{enumerate}
Likewise, the condition that every $p$th root of $1_G$ is trivial is also equivalent to two other conditions:
\begin{enumerate}
    \item for all $k$, every $p^k$th  root of $1_G$ is trivial;
    \item for all $k>0$, every $a \in G$ has exactly one $p^k$th  root.
\end{enumerate}
\end{rem}
\noindent From Remark~\ref{Rem: the remark on roots}, we easily deduce the following:
\begin{rem} \label{rem: Up is a model of GM}
 For every $p$, $\uu_{(p)}$ is a model of $\thgp$, and so is $\ff^\times$ when $\text{char}(\ff)=p$. Moreover, if $G$ is a model of $\thgp$, then the group of torsion elements of $G$ is isomorphic to $\uu_{(p)}$. The group $\uu$ is a model of $\thgo$ and is isomorphic to the group of torsion elements of any $\thgo$-model.
\end{rem}

\noindent Lemma~\ref{lem: property multiplicative group} confirms that our candidate $\thgg$ at least  meets the basic requirements:

\begin{lem} \label{lem: property multiplicative group}
If $G$ is the multiplicative group of a model of $\ACF$, then $G \models \thgg$.
Similar statements hold for $\ACF_p$ together with $\thgp$ for an arbitrary $p$ and for $\ACF_0$ together with $\thgo$.
\end{lem}
\begin{proof}
It is easy to see that if $G$ is the multiplicative group of a prime model of $\ACF$, then conditions (G0), (G1), and (G2) are satisfied. Hence, the first statement follows from the fact that $\ACF$ is model complete. The proof of the second statement is similar.
\end{proof}

\noindent Suppose $B$ is a subset of $G$, and $t(x)$ and $t'(x)$ are $\lag(B)$-terms. Then we call the atomic formula $t(x) = t'(x) $  a {\bf multiplicative equation over $B$}. A multiplicative equation over $B$ is {\bf trivial} if it defines in every abelian group $G'$ extending $\langle B \rangle$ the set $(G')^m$.  If $a \in G^m$ does not satisfy any nontrivial multiplicative equation over $B$, we say that $a$ is {\bf multiplicatively independent} over $B$.

\meno  Proposition~\ref{prop: GMproperties} below is the  ``almost'' quantifier elimination result we promised. This can be seen as folklore and can be obtained as a consequence of the characterization of elementary embeddings of abelian groups \cite{Eklof} and the quantifier reduction for abelian groups \cite[p. 46]{MikePrest}. Since the situation is relatively simple, we briefly indicate a direct proof:

\begin{prop} \label{prop: GMproperties}
 For each $p$, the theory $\thgp$ is complete and admits quantifier elimination.  A similar statement holds for $\thgo$. However, $\mathrm{GM}$ is not model complete.
\end{prop}
\begin{proof}

 We will only prove the first statement for $\thgp$ with $p$ fixed as the proof for $\thgo$ is very similar. By Remark~\ref{rem: Up is a model of GM}, $\uu_{(p)}$ is an $\lag$-substructure of every model of $\thgp$, so completeness will follow from quantifier elimination. 
Using a standard test for quantifier elimination, we need to show the following: if $G$ and $ G'$ are models of $\thgp$ such that $G'$ is $|G|^+$-saturated, and $f$ is a partial $\lag$-isomorphism from $G$ to $G'$ (i.e., $f$ is an $\lag$-isomorphism from an $\lag$-substructure of $G$ to an $\lag$-substructure of $G'$) such that $\text{Domain}(f) \neq G$, then there is a partial $\lag$-isomorphism from $G$ to $G'$ which properly extends $f$. 

 In each of the following cases, we will obtain $a$ in $G \setminus \text{Domain}(f)$ and $a'$ in $G' \setminus \text{Image}(f)$. A proper extension of $f$ can then be defined by
$$a^kb \mapsto (a')^k f(b) \ \text{ for } k \in \zz \text{ and } b \in \text{Domain}(f).$$ We will leave the reader to check that the function is well-defined and is a partial $\lag$-isomorphism from $G$ to $G'$.

Suppose $l$ is a prime number, and $a \in G \setminus \text{Domain}(f)$ is a nontrivial $l$th root of $1_G$. As $G$ satisfies (GC$_p$), $l \neq p$. Since $G$ and $G'$ both satisfy (G0), $\text{Domain}(f)$ and $\text{Image}(f)$ contain no nontrivial $l$th roots of $1_G$ and $1_{G'}$ respectively. We can then choose $a' \in G' \setminus \text{Image}(f)$ to be an $l$th root of $1_{G'}$, which must exist because $G'$ satisfies (GC$_p$) and (G2). 

Now suppose $\text{Domain}(f)$ contains all roots of $1_G$ with prime order,  $l$ is a prime and $a \in G \setminus \text{Domain}(f)$ is such that $a^l \in \text{Domain}(f)$. If $b$ is another $l$th root of $a^l$, then $ab^{-1}$ is an $l$th root of $1_G$. Hence, $\text{Domain}(f)$ contains no $l$th root of $a^l$, and  $\text{Image}(f)$ contains no $l$th root of $f(a^l)$. We then choose $a'$ to be an $l$th root of $f(a^l)$ which must exist because $G'$ satisfies (G1). 

The last case is when $\text{Domain}(f)$ is divisibly closed in $G$, and $a \in G \setminus \text{Domain}(f)$. Using the fact that $G'$ is $|G|^+$-saturated, we obtain  $a' \in G'$ which is multiplicatively independent over $\text{Image}(f)$.

For the last statement, note that  both $\uu_{(p)}$ and $\uu$ are models of $\thgg$, and $\uu_{(p)}$ is a substructure of $\uu$, but $\uu_{(p)}$ is not  not an elementary substructure of $\uu$.
\end{proof}
\noindent Fact~\ref{Fact: Circularly orderable groups} is an easy consequence of \v{Z}eleva's characterization of circularly orderable groups\cite{Zeleva} and Levi's characterization of linearly orderable abelian group \cite{Levi}:
\begin{fact} \label{Fact: Circularly orderable groups}
An abelian group is circularly orderable if and only if it satisfies (G0). 
\end{fact}

\noindent Combining Proposition~\ref{prop: GMproperties} and Fact~\ref{Fact: Circularly orderable groups} confirms the validity of our candidate $\thgg$:
\begin{cor}
Every model of $\thgg$ is elementarily equivalent to both the multiplicative group of an algebraically closed field and the multiplicative group of a model of $\ACFO^-$.
\end{cor}
\noindent Many other model-theoretic properties of the theory $\thgg$ are also immediate:
\begin{cor} \label{Cor: strongly minimal}
The theory $\thgg$ is strongly minimal.
\end{cor}
\noindent Hence, definable sets, types, and elements in a model of $\thgg$ can be given a canonical dimension $\text{mdim}$ which coincides with Morley rank and the $\acl_\text{m}$-dimension; see \cite{Marker} for details. Proposition~\ref{prop: GMproperties} also yields:
\begin{cor} \label{dimensionandmultiplicativelarge}
Suppose $G$ is a model of  $\thgg$,  $B$ is a subset of $G$, and $a$ is in $G^m$. Then $\mathrm{mdim}(a|B) < m$ if and only if  $a$ is multiplicatively dependent over $B$.
 \end{cor}

\subsection{Circularly ordered multiplicative groups of $\ACFO^-$-models} \label{Circularly ordered multiplicative groups of ACF-models}
We next consider the theory of  circularly ordered multiplicative groups of models of $\ACFO^-$. We want to show that that this theory ``almost'' has a model companion and obtain an axiomatization for this model companion along the way.

\meno Throughout Section~\ref{Circularly ordered multiplicative groups of ACF-models}, we adopt the notational conventions of Section~\ref{Sec: Multiplicative groups of algebraically closed field}. Moreover, $G$ is assumed to be circularly orderable, and $(G, \triangleleft)$ ranges over the circularly ordered multiplicative abelian groups.
For each $(G, \triangleleft)$, we define the linear order $\lessdot$ on $(G, \triangleleft)$ by setting $1_G \lessdot a$ for all $a \in G\setminus\{1_G\}$ and 
$$ a \lessdot b \text{ if and only if } \triangleleft(1_G,a ,b) \quad \text{ for } a, b \in G \setminus \{ 1_G\}. $$
When $G$ is $\mathbb{T}$, $\uu$, or $\uu_{(p)}$, we let  $\triangleleft$ denotes the clockwise circular orders on the respective sets. From Fact~\ref{Fact: Circularly orderable groups}, we can easily deduce the following:
\begin{rem} 
For given $(G, \triangleleft)$  and finite subgroup $A$ of $G$, if $a \in A \setminus\{1_G\}$ is minimal with respect to $\lessdot$, then $A = \la a \ra$.
\end{rem}

\noindent Let $\lam = \lag \cup \{ \triangleleft\}$ be the language of circularly ordered abelian groups. Let $\thmm^-$ be the theory whose models are $(G, \triangleleft)$ such that $G \models \thgg$, or equivalently, $G$ satisfies (G1) and (G2) (as (G0) is automatic by Fact~\ref{Fact: Circularly orderable groups}). Let $\thmp^- = \thmm^- \cup \thgp$ for all $p$, and let  $\thmo^- = \thmm^- \cup \thgo$.  We show below that $\thmm^-$ is an axiomatization of the theory of circularly ordered multiplicative groups of algebraically closed fields: 

\begin{lem}
An $\lam$-structure $(G, \triangleleft)$ is a model of $\thmm^-$ if and only if $(G, \triangleleft)$ is elementarily equivalent to the circularly ordered group of an $\ACF$-model. Similar statements hold for $\thmp^-$ together with $\ACF_p$ for an arbitrary $p$ and $\thmo^-$ together with $\ACF_0$.
\end{lem}
\begin{proof}
The backward implication of the first statement follows immediately from Lemma~\ref{lem: property multiplicative group}. For the forward implication of the first statement, suppose  $(G, \triangleleft)$ is a model of $\thmm^-$. We assume further that $(G, \triangleleft)$ is a model of $\thmp^-$ and omit the proof of the similar case where $(G, \triangleleft)$ is a model of $\thmo^-$.
 Replacing $(G, \triangleleft)$ by an elementary extension if necessary, we can arrange that $|G| =\kappa > \aleph_0$. By Corollary~\ref{Cor: strongly minimal}, $\thmp^-$ is $\kappa$-categorical. Hence, $G$ is is isomorphic to the multiplicative group $G'$ of a model of $\ACF_p$ of size $\kappa$. Pushing forward $\triangleleft$ by the isomorphism we obtain a circular orderding  $\triangleleft'$ on $G'$ such that  $(G', \triangleleft')$ is $\lam$-isomorphic to $(G, \triangleleft)$.  This also proved the second statement.
\end{proof}

\noindent A rather awkward aspect dealing with $(G, \triangleleft)$ comes from the fact that $\lessdot$ is not invariant under translation. We will consider here a partial rectification.
The {\bf winding number} $W(a_1, \ldots, a_n)$  of  $(a_1, \ldots, a_n) \in G^n$ is defined to be the cardinality of the set
$$  \big\{ k \mid  1 \leq k \leq n-1, \prod_{i=1}^{k+1} a_i \lessdot \prod_{i=1}^{k} a_i \big\}. $$
It is intuitively the number of times the sequence $a_1, a_1a_2, \ldots, \prod_{i=1}^{n-1} a_i, \prod_{i=1}^n a_i$ ``winds around the circle''. 
If $a_1 =\ldots = a_n =a$, we also denote $W(a_1, \ldots, a_n)$ as $W_n(a)$.

\begin{rem} \label{rem: Winding number rectify noninvariant}
Suppose $a$ and $b$ in $G$ satisfy $a \lessdot b$. Then for all $c \in G$, either $ac \lessdot bc$ or $W(a,c) <W(b,c)$. In this sense, the notion of winding number accounts for the noninvariance of $\lessdot$.  
\end{rem}

\noindent For $a \in G$, we say that $a$ is {\bf $n$-divisible with winding number $r$} if $a$ has an $n$th root $b$ with $W_n(b)=r$. Remark~\ref{rem: Winding number rectify noninvariant} then gives us: 

\begin{rem} \label{r-divisibility is unique}
For any given $(G, \triangleleft)$, every element $a \in G$ has at most one $n$th root $b$ such that $W_n(b)=r$. So if there are distinct $b_1, \ldots, b_n$ such that $b_i^n =a$ for all $i \in \{1, \ldots, n\}$, then for each $r \in \{1, \ldots, n\}$ there is exactly one $i \in \{1, \ldots, n\}$ such that $W_n(b_i)=r$. 
\end{rem}

\noindent  Let $\thmm$ be the $\lam$-theory such that its models are $(G, \triangleleft)$ with $G \models \thgg$ and the following density condition is satisfied:
\begin{enumerate}
    \item[(GD)] for any given $n$, $r \in \{0, \ldots, n-1\}$, and $c$ and $d$ in $G$, there is $a \in G$ such that $\triangleleft(c,a,d)$ and $a$ is $n$-divisible with winding number $r$.
\end{enumerate}
The theory $\thmm$ is our candidate for the ``almost'' model companion of $\thmm^-$. Also set $\thmp = \thmm \cup \thmp^-$ for an arbitrary $p$, and set $\thmo = \thmm \cup \thmo^-$.

\meno To handle circularly ordered groups, it is convenient to ``linearize'' them; see also \cite{Ayhan} and \cite{Giraudet} for related material.  Let $(H, <)$ be a linearly ordered additive group with identity element $0_H$,  and  let $\omega \in H$ be a distinguished positive element  such that $( n\omega)_{n>0}$ is cofinal in $(H, <)$ (i.e., for every $\alpha \in H$, $\alpha < n\omega$ for sufficiently large $n$). For every $k$, set 
 $$[k, k+1)_H = \{ \alpha  \in H
\mid  k\omega \leq \alpha < (k+1) \omega\}. $$
\noindent A surjective group homomorphism $e : H \to G$ with kernel $\la \omega \ra$ is a
{\bf covering map}  from $(H, \omega, <)$ to $(G, \triangleleft)$ if for all $n$ and all $\alpha, \beta, \gamma \in [n, n+1)_H$, we have that $\triangleleft(e(\alpha), e(\beta), e(\gamma))$ is equivalent to
$$  \alpha<\beta<\gamma \quad \text{or} \quad \beta< \gamma< \alpha \quad \text{or}\quad  \gamma< \alpha < \beta.$$
If there is a covering map from $(H, \omega,<)$ to $(G,\triangleleft)$, we call $(H, \omega,<)$ a {\bf universal cover} of $(G,\triangleleft)$.

\begin{rem} \label{Rem: Cover backward}
Suppose $(H, \omega,<)$ is  as described in the preceding paragraph. Then the above definition also allow us to construct $(G, \triangleleft)$ such that  $(H, \omega,<)$ is a universal cover of $(G, \triangleleft)$.
\end{rem}

\noindent The examples in the  following remark will hopefully make this notion concrete:

\begin{rem}\label{Rem: Zp, qq, rr} 
Let the additive groups $\rr$, $\qq$, and $\zz_{(p)}$ be equipped with their natural orders $<$. With $\alpha \mapsto e^{2 \pi i \alpha}$ the covering map, we have the following:
\begin{enumerate}
    \item $(\rr, 1, <)$ is a universal cover of $(\mathbb{T}, \triangleleft)$; 
    \item $(\qq, 1, <)$ is a universal cover of $(\uu, \triangleleft)$; 
    \item $(\zz_{(p)}, 1, <)$ is a universal cover of $(\uu_{(p)}, \triangleleft)$.
\end{enumerate}
\end{rem}

\noindent The lemma below illustrates the advantage of having a universal cover.

\begin{lem} \label{Universalcoveradvantage}
Suppose $(H, \omega, <)$ is a universal cover of $(G, \triangleleft)$ with $e$ the covering map, $\alpha_1,\ldots, \alpha_n$ are in $[0,1)_H$, and $a_i = e(\alpha_i)$ for $i \in \{1, \ldots, n\}$. Then
$$ W(a_1, \ldots, a_n) =r \text{ if and only if } \alpha_1 + \ldots+ \alpha_n \in [r, r+1)_H. $$
\end{lem}
\begin{proof}
It follows from the definition of a universal cover that $\sum_{i=1}^k \alpha_i \in [l, l+1)_H$ and $\sum_{i=1}^{k+1} \alpha_i \in [l+1, l+2)_H$  if and only if $\prod_{i=1}^{k+1} a_i \lessdot \prod_{i=1}^k a_{i}$. The desired conclusion follows.
\end{proof}

\noindent Applying Lemma~\ref{Universalcoveradvantage} into the setting where $a_1 =\ldots =a_n =a$ yields:

\begin{cor} \label{Windinganddivisibility}
Suppose $(H, \omega, <)$ is a universal cover of $(G, \triangleleft)$ with $e$ the covering map, $a$ is in $G$, $n\geq 1$,  $r$ is in $\{0, \ldots, n-1\}$, and $\alpha \in [r, r+1)_H$ is such that $e(\alpha)  =a$. Then the following are equivalent:
\begin{enumerate}
    \item[\ir] $a$ is $n$-divisible with winding number $r$;
    \item[\iir] $\alpha$ is $n$-divisible.
\end{enumerate}
\end{cor}

\noindent  We can view such $(H, \omega, <)$ as a structure in a language $\laa$ consisting of function symbols for $0$, $\omega$, and $+$ and a relation symbol for $<$. 
It turns out that the convenience of a universal cover is something we can always afford. Moreover, we get it  partially definably:

\begin{lem}\label{prop:cover}
Every $(G, \triangleleft)$ has a universal cover $(H, \omega, <)$. Moreover, there is an $\lam$-isomorphic copy $(\tilde{G}, \tilde{\triangleleft})$ of  $(G, \triangleleft)$ such that the underlying set of $\tilde{G}$ is $[0, 1)_H$, and the multiplication on $\tilde{G}$ and $\tilde{\triangleleft}$ can be defined by $\laa$-formulas whose choice is independent of the choice of $(G, \triangleleft)$ and the choice of $(H, \omega, <)$.
\end{lem}
\begin{proof}
Set $H = \zz \times G$, and define 
$$(k, a)+ (k', a') = (k+k'+ W(a, a'), aa')
 $$
 for $(k, a)$ and $(k', a')$ in $H$. 
Let $<$ be the lexicographic product of the usual order on $\zz$ and the linear order $\lessdot$ on $G$. Set $0_H = (0_\zz, 1_G)$ and $\omega = (1_\zz, 1_G)$. We can easily check that $(H,\omega,+,<)$ is a universal cover of $G$.
For $a,a' \in [0_H,\omega)_H$, set
$$ a \ \tilde{\times} \  a'=
\begin{cases}
a + a' & \text{ if } a + a' \in [0,1)_H ,\\
a + a' - \omega & \text{ otherwise}.
\end{cases}
$$
Define $\tilde{\triangleleft}$ by setting $\tilde{\triangleleft}(a,b,c)$ for any $a,b,c \in [0,1)_H$ such that $a < b < c$ or $b < c < a$ or $c < a < b$.
It is easy to see the quotient map $H \to G$ induces an isomorphism from  $([0,1)_H,\tilde{\times},\tilde{\triangleleft})$ to $(G,\triangleleft)$.
\end{proof}

\noindent The universal cover notion is functorial in the following sense:

\begin{lem} \label{lem:  functoriality of covering}
Suppose $(H, \omega, <)$ is a universal cover of $(G, \triangleleft)$ with convering map $e$, and $(H', \omega', <')$ is a universal cover of $(G', \triangleleft')$ with convering map $e'$. Then we have the following:
\begin{enumerate}
    \item[\ir] if  $g$ is an $\lam$-embedding from $(G, \triangleleft)$ to  $(G', \triangleleft')$, then there is a unique $\laa$-embedding $h$ from $(H, \omega, <)$ to $(H', \omega', <')$ such that the diagram below commutes:
    \begin{center}\begin{tikzcd}
    H'\ar[r,"e'"]   &  G'\\
    H\ar[u,"h"']\ar[r,"e"]&G\ar[u,"g"']
    \end{tikzcd};\end{center}
        in particular, $e$ is the unique covering map from $(H, \omega, <)$ to $(G, \triangleleft)$, and any two universal coverings of $(G, \triangleleft)$ are isomorphic as $\laa$-structures;
    \item[\iir] if  $h$ is an $\laa$-embedding from $(H, \omega, <)$ to  $(H', \omega', <')$, then there is a unique $\lam$-embedding $g$ from $(G, \triangleleft)$ to $(G', \triangleleft')$ such that the same diagram above commutes.
\end{enumerate}
\end{lem}

\begin{proof}
For (i), let $h: H \to H'$ be such that  $ \alpha \in [k, k+1)_H$ is mapped to the unique $\beta \in [k, k+1)_{H'}$ with $g \circ e(\alpha) = e'(\beta)$. For (ii), let $g: G \to G'$ be such that  $e(\alpha)$ is mapped to $e'\circ h(\alpha)$ for $\alpha \in H$. It is easy to check that $h$  and $g$ are as desired.
\end{proof}

\noindent We extend the ``linearization'' procedure to theories $\thmm^-$ and $\thmm$. Let $\thaa^-$ be an $\laa$-theory such that an $\laa$-structure $(H, \omega, <)$ is a model of $\thaa^-$ if and only if
$(H, <)$ is a linearly ordered additive abelian group, $\omega$ is a positive element in $H$, and the the following additional two properties are satisfied:
\begin{enumerate}
    \item[(H1)] for each $n$ and $\alpha \in H$, there is at least one $r \in \{0, \ldots, n-1\}$ such that $\alpha+ r\omega$ is $n$-divisible;
       \item[(H2)]  for any prime numbers $p$ and $l$, $\omega$ is either $p$-divisible or $l$-divisible.
\end{enumerate}
Note that (H1) and (H2) correspond to (G1) and (G2). There is no (H0) because (G0) is trivial in our current setting. The condition that $\omega$ is cofinal in $H$ cannot be included here as it is not first-order.
For a given $p$, let $\thap^-$ be the $\laa$-theory  whose models are the $(H, \omega, <) \models \thaa^-$ which satisfy the addition condition:
\begin{enumerate}
    \item[(HC$_{p}$)] $\omega$ is not $p$-divisible. 
\end{enumerate}
 We also let let $\thao^-$ be the $\laa$-theory  whose models are the $(H, \omega, <) \models \thaa^-$ which satisfy the additional condition:
\begin{enumerate}
    \item[(HC$_{0}$)]  for all prime numbers $l$, $\omega$ is $l$-divisible.
\end{enumerate}
Let $\thaa$ be an $\laa$-theory whose  models are  the $(H, \omega, <) \models \thaa^-$ which also satisfy the additional condition: 
\begin{enumerate}
    \item[(HD)] for any given $n$ and $\beta, \gamma \in H$ with $\beta < \gamma$, there is $\alpha \in H$ such that $\alpha$ is $n$-divisible and   $\beta < \alpha < \gamma$.
\end{enumerate}
Finally, set $\thap = \thaa \cup \thap^-$ for each $p$, and $\thao = \thaa \cup \thao^-$; in fact, $\thao$ is just the theory of divisible ordered abelian groups. The next Lemma explains precisely what it means by saying that these are ``linearization'' of $\thmm^-$ and $\thmm$:

\begin{lem}\label{lem:GMOtoHAO}
Suppose $(H, \omega, <)$ is a universal cover of $(G, \triangleleft)$. Then we have:
\begin{enumerate}
    \item[\ir] for all $p$, $(H, \omega, <) \models \thaa^-$ if and only if  $(G, \triangleleft) \models \thmm^-$. Similar statements hold for $\thap^-$ together with $\thmp^-$ and $\thao^-$ together with $\thmo^-$;
    \item[\iir] for all $p$, $(H, \omega, <) \models \thaa$ if and only if  $(G, \triangleleft) \models \thmm$. Similar statements hold for $\thap$ together with $\thmp$ and $\thao$ together with $\thmo$.

\end{enumerate}
\end{lem}

\begin{proof}
All these statements are immediate consequences of Corollary~\ref{Windinganddivisibility}.
\end{proof}

\noindent Lemma~\ref{lem:GMOtoHAO} allows us to deduce results for $\thmm^-$-models and $\thmm$-models from generally much easier results for $\thaa^-$-models and $\thaa$-models. Below is the first demonstration of its usefulness: 

\begin{lem} \label{Zp is model of HAO-}
Let  $(\zz_{(p)}, 1, <)$ and $(\qq, 1, <)$ be as in Remark~\ref{Rem: Zp, qq, rr}. Then we have $(\zz_{(p)}, 1, <) \models \thap$ and $(\qq, 1, <) \models \thao$. Moreover, there is a unique  $\laa$-embedding of $(\zz_{(p)}, 1, <)$ into every $\thap^-$-model and a unique $\laa$-embedding of $(\qq, 1, <)$ into every $\thao^-$-model.
\end{lem}
\begin{proof}
 It is easy to verify that $(\zz_{(p)}, 1, <)$ is a model of $\thap^-$. Since $\zz_{(p)}$ is dense in $\rr$ with respect to the natural order, it follows that  $(\zz_{(p)}, 1, <)$ is a model of $\thap$. Suppose $(H, \omega, < )$ is a model of $\thap^-$. Then the subgroup of $H$ generated by $\omega$ is an isomorphic copy of $\zz_{(p)}$. This gives us an $\laa$-embedding of $(\zz_{(p)}, 1, <)$ into $(H, \omega, < )$. This $\laa$-embedding is unique as any such $\laa$-embedding must send $1$ to $\omega$. The statements for $(\qq, 1, <)$ can be proven similarly.
 \end{proof}
\noindent Combining with Lemma~\ref{lem:  functoriality of covering} and Lemma~\ref{lem:GMOtoHAO}, we get:
\begin{prop} \label{Up is model of GMO-}
We have $(\uu_{(p)}, \triangleleft) \models \thmp$ and $(\uu, \triangleleft) \models \thmo$. Moreover, there is a unique  $\lam$-embedding of $(\uu_{(p)}, \triangleleft)$ into every   $\thmp^-$-model and a unique $\lam$-embedding of $(\uu, \triangleleft)$ into every $\thmo^-$-model.
\end{prop}

\noindent Suppose $\sigma$ is an $\lag$-automorphism of $\uu_{(p)}$. Define  $\triangleleft_\sigma$ to be the image of the clockwise circular order $\triangleleft$ under $\sigma$. From Lemma~\ref{Up is model of GMO-}, we have the following:

\begin{cor} \label{Up is model of GMO-2}
Every circular order on $\uu_{(p)}$ is equal to $\triangleleft_\sigma$ for a unique $\lag$-automorphism $\sigma$ of $\uu_{(p)}$.
\end{cor}

\noindent For   an injective group homomorphism $\chi: \ff^\times \to \mathbb{T}$, define the circular order $\triangleleft_\chi$ to be the pullback of $\triangleleft$ via $\chi$. Note that  $\chi(\ff^\times) = \uu_{(p)}$ as a subgroup of $\bbT$. So applying Corrollary~\ref{Up is model of GMO-2}, we deduce:

\begin{cor} \label{Up is model of GMO-3}
Every multiplicative circular order on $\ff$ is equal to $\triangleleft_\chi$ for a unique  injective group homomorphism $\chi: \ff^\times \to \mathbb{T}$.
\end{cor}

\noindent Toward showing that $\thmm$ is ``almost'' the model companion of $\thmm^-$, we first show the ``linearized'' version of the result. This is also folklore \cite{Weis}, but not everything we want is written down, so we briefly indicate a proof.

\begin{lem} \label{prop: HAOproperties}
For each $p$, the theory $\thap$ is complete and is the model companion of $\thap^-$. A similar statement holds for $\thao$ and $\thao^-$. The theory $\thaa$ is not model complete.
\end{lem}

\begin{proof}

To show the first statement, we require some preparation. For each model $(H, \omega, <)$ of $ \thap^-$, define the family $D=(D_n)_{n\geq 1}$ of unary relations  by setting
$$D_n = \{ \alpha \in H \mid \text{ there is } \beta \in H \text { such that } n\beta =\alpha.\}$$ 
 Then such $(H,\omega, <, D)$  is naturally a structure in a language $\laa^\diamondsuit$ extending $\laa$ by adding a family of unary relation symbols for  $D$. The theory $\thap^-$ and $\thap$ can be naturally expanded to $\laa^\diamondsuit$-theories by adding the obvious axioms defining such $D$. Note that when $(H, \omega, <) \models \thap^-$ and $D=(D_n)_{n>0}$ are as above, we also have $$D_n = \{ \alpha \in H \mid \text{ for all } \beta \in H,  \bigwedge_{r=1}^{n-1} n\beta \neq \alpha+r\omega.\}$$
 It follows that such $D_n$ is both universally and existentially definable. Moreover, we can choose the formula defining such $D_n$ independent of the choice of  $(H, \omega, <)$. Thus, the problem is reduced to showing that the natural $\laa^\diamondsuit$-expansion of $\thap$ is complete, admits quantifier elimination,  and is the model companion of the natural $\laa^\diamondsuit$-expansion of $\thap^-$.

It follows from Lemma~\ref{Zp is model of HAO-} that $(\zz_{(p)}, 1, <, D)$ can be canonically viewed as a $\laa^\diamondsuit$-substructure of any model of the natural $\laa^\diamondsuit$-expansion of $\thap$. Hence,it suffices to show the following: if 
$(H, \omega, <, D)$ is the natural $\laa^\diamondsuit$-expansion of a model of $\thap^-$, $(H', \omega', <', D')$ is the natural $\laa^\diamondsuit$-expansion of a model of $\thap$ and is moreover $|H|^+$-saturated, and $f: H \to H'$ is a partial $\laa^\diamondsuit$-isomorphism from $(H, \omega, <, D)$ to $(H', \omega', <', D')$ with $\text{Domain}(f) \neq H$, then we can find a partial $\laa^\diamondsuit$-embedding which properly extends $f$. 

It is easy to reduce to the case where $\text{Domain}(f)$ is a divisibly closed subgroup of $H$. Let $\alpha \in  H \setminus \text{Domain}(f)$. If $\alpha -r\omega$ is $p^k$-divisible and $\beta < \alpha < \beta'$ for $\beta$ and $\beta'$ in $\text{Domain}(f)$, then we can find  $\alpha'$ in $H' \setminus \text{Image}(f)$  such that $\alpha' -r\omega'$ is $p^k$-divisible and $f(\beta) < \alpha' < f(\beta')$ using the fact that $(H', \omega', <')$ satisfies (HD). As $H'$ is $|H|^+$-saturated, we can arrange that $\alpha'$ satisfies all such conditions simultaneously. Let $g$ be the obvious extension of $f$ sending $\alpha$ to $\alpha'$. It is easy to check that $g$ is as desired.

The proof of the second statement is similar to the proof of the first statement.
Note that $(\zz_{(p)}, 1, <)$ is an $\laa$-substructure of $(\qq, 1, <)$, and both are models of $\thaa$, but the former is not an elementary substructure of the latter. So $\thaa$ is not model complete.
\end{proof}

\begin{prop} \label{prop: GMOproperties}
For each $p$, the theory $\thmp$ is complete and is the model companion of $\thmp^-$. A similar statement holds for $\thmo$ and $\thmo^-$. However, $\thmm$ is not model complete.
\end{prop}

\begin{proof}

By Proposition~\ref{Up is model of GMO-},  $(\uu_{(p)}, \triangleleft)$ is a model of $\thmp$ and is an $\lam$-substructure of any model of $\thmp$. Hence to get the completeness of $\thmp$, it suffices to show that $\thmp$ is model complete. Suppose $(G, \triangleleft)$ and $(G', \triangleleft)$ are models of $\thmp$ and that the former is a substructure of the latter. By Lemma~\ref{lem:GMOtoHAO}, the universal covers $(H, \omega, <)$ and  $(H', \omega, <)$ of $(G, \triangleleft)$ and $(G', \triangleleft)$ are models of  $\thap$. It follows from Lemma~\ref{lem:  functoriality of covering} and Lemma~\ref{prop: HAOproperties} that $(H, \omega, <)$ can be viewed as an elementary substructure of  $(H', \omega, <)$. Combining this with the second part of Lemma~\ref{prop:cover}, we learn that  $(G,  \triangleleft)$ is an elementary substructure of $(G', \triangleleft)$.

Next we show that every model of $\thmp^-$ can be embedded into a model of $\thmp$. Suppose $(G, \triangleleft)$ is a model of $\thmp^-$. By Lemma~\ref{lem:GMOtoHAO}, the universal cover $(H, \omega, <)$ of  $(G, \triangleleft)$ is a model of $\thap^-$. Hence, it follows from Lemma~\ref{prop: HAOproperties} that  $(H, \omega, <)$ has an extension $(H', \omega,  <)$ which is a model of $\thap$. Construct $(G', \triangleleft)$ as mentioned in Remark~\ref{Rem: Cover backward}. Then $(G', \triangleleft)$ is a model of $\thmp$ by Lemma~\ref{lem:GMOtoHAO} and $(G', \triangleleft)$ can be considered a substructure of $(G', \triangleleft)$ by Lemma~\ref{lem:  functoriality of covering}.

The second statement can be proved similarly. The third statement can be deduced from Lemma~\ref{prop: HAOproperties} using similar ideas. It can also be observed directly by looking at $(\uu_{(p)}, \triangleleft)$ and $(\qq, \triangleleft)$.
\end{proof}

\begin{rem}
The theory $\thao$ is just the theory of divisible ordered abelian groups, so $\thao$ has quantifier elimination. For the second statement of Proposition~\ref{prop: GMOproperties}, we can also have that  $\thmo$ admits quantifier elimination.  Since we will not use this later on, we leave it to the interested reader.
\end{rem}

\section{Model companion of $\ACFO^-$} \label{sec: Model  companion of ACFO-}

\noindent We will establish that  $\ACFO^-$ has a model companion in two steps:
\begin{enumerate}
      \item obtaining a characterization of the existentially closed models  of $\ACFO^-$ following the ideas in    \cite{ChatPillay};
    \item showing that the  class of $\ACFO^-$-models satisfying the characterization in (1) is elementary by using model-theoretic/geometric properties of the reducts of $\ACFO^-$ to the language of rings and the language of circularly ordered multiplicative groups.  
\end{enumerate}

\subsection{Geometric characterization of the existentially closed models} Intuitively, in an existentially closed model $(F, \triangleleft)$ of $\ACFO^-$, the field $F$ interacts ``randomly'' with the circularly ordered abelian group $(F^\times, \triangleleft) \models \thmm$ over their ``common reduct'' $F^\times$. In this section, we will make precise this intuition through a ``geometric characterization'' and then verify its correctness.

\meno We keep the the notational conventions of Section~\ref{Sec: Multiplicative groups of algebraically closed field} and Section~\ref{Circularly ordered multiplicative groups of ACF-models}.
Suppose $(G, \square)$ is an  $L$-structure expanding $G$. For convenience, we call a set $X \subseteq G^m$ which is defined in $(G, \square)$ by a quantifier-free $L(G)$-formula a {\bf qf-set} in $(G, \square)$. For $X \subseteq G^m$ definable in $(G, \square)$ and an elementary extension 
$(G',\square)$ of $(G, \square)$, let $X(G') \subseteq (G')^m$ be the set defined in $(G', \square)$ by an arbitrary $L(G)$-formula  $\phi(x)$ that defines $X$.

\meno We first correct a minor issue: the group $F^\times$ is, strictly speaking, not a reduct of $F$,  as $0$ is not an element of $F^\times$.  Set
$$ \Sigma_{n+1} = \{ (a_1, \ldots, a_{n+1}) \in (F^\times)^{n+1} \mid a_1+ \ldots+a_n =a_{n+1} \}$$
and let $\Sigma = (\Sigma_{n+1})$. We call $(F^\times, \Sigma)$ the {\bf punctured field} associated to $F$. Then $(F^\times, \Sigma)$ is naturally a structure in the language $\lar^\times = \lag \cup \{ \Sigma_{n+1}\}$. The group $F^\times$ is now an honest reduct of $(F^\times, \Sigma)$.

\meno We will see in the proof of Lemma~\ref{genericimpliesexistentiallyclosed} a more  substantial advantage working with $\lar^\times$ instead of $\lar$, namely,  $\lar^\times$ expands $\lag$ only by relation symbols and not by function symbols.

\begin{rem} \label{rem: adding0}
The following  ``adding 0'' procedure allows us to recover an isomorphic copy of a field from its associated punctured field, but the procedure is applicable to any $\lar^\times$-structure expanding a multiplicative abelian group. Starting with an $\lar^\times$-structure $(G, \Sigma)$, set $F = G \cup \{0\}$,  define $+$ on $F^2$ by pretending that $G$ is $F^\times$  (i.e., $0+0=0$, $a+0=0+a =a$ for $a \in G$, $a+b = c$ for $a$ and $b$  in $G$ if $c$ is the unique element of $G$ satisfying  $\Sigma_3(a,b,c)$, and $a+b=0$ for the remaining cases), and define $\times$ on $F^2$ similarly.
\end{rem}
 
\noindent   As immediate consequence of  Remark~\ref{rem: adding0}, we get:
\begin{rem} \label{rem: definable sets in a field and its associated punctured field}
Suppose $F$ is a field and $(F^\times, \Sigma)$ is its associated punctured field. Then $X \subseteq (F^\times)^m$ is definable in $F$ if and only if $X$ is definable in $(F^\times, \Sigma)$.
\end{rem}

\noindent From Remark~\ref{rem: adding0}, it is also easy to find an $\lar^\times$-theory whose models are precisely the punctured fields. Likewise, we get $\lar^\times$-theories $\ACF^\times$, $\ACF^\times_p$ for every $p$, and  $\ACF^\times_0$ whose models are punctured models of  $\ACF^\times$, $\ACF^\times_p$, and $\ACF^\times_0$ respectively. The basic model theory $\ACF^\times$ can be obtained:

\begin{lem}\label{Lem: Theory of punctured field}
The theory $\ACF^\times$ admits quantifier elimination and is the model companion of the theory of punctured fields. The theories $\ACF^\times_p$ for various $p$ and $\ACF^\times_0$ are the only completions of $\ACF^\times$.
\end{lem}

\begin{proof}
These statements are easy consequences of Remark~\ref{rem: adding0}, Remark~\ref{rem: definable sets in a field and its associated punctured field}, the quantifier elimination of $\ACF$, and the fact that  $\ACF_p$ for various $p$ and $\ACF_0$ are the only completions of $\ACF$.
\end{proof}

\noindent Let  $(F, \triangleleft)$ be an $\ACFO^-$-model and $(F^\times, \Sigma)$ the punctured field associated to $F$. We call $(F^\times, \Sigma, \triangleleft)$ the  {\bf punctured $\ACFO^-$-model} associated to $(F, \triangleleft)$.  Then $(F^\times, \Sigma, \triangleleft)$ is a structure in the language $\lac^\times = \lar^\times \cup \lam$. We define punctured $\ACFO^-_p$-models for various $p$ and punctured $\ACFO^-_0$-models likewise.

\meno By the  discussion on ``adding 0'' in Remark~\ref{rem: adding0}, it is easy to see that there is an $\lac^\times$-theory whose models are precisely the punctured $\ACFO^-$-models. 
We say that a punctured $\ACFO^-$-model is {\bf existentially closed} if it is an existentially closed model of this theory. 

\meno The following Lemma allows us to trade existentially closed $\ACFO^-$-models with existentially closed punctured $\ACFO^-$-models:

\begin{lem} \label{Existentially closed models exchange}
An $\ACFO^-$-model is existentially closed if and only if its associated punctured $\ACFO^-$-model is existentially closed. \end{lem}
\begin{proof}
Let $(F^\times, \Sigma, \triangleleft)$ be the punctured $\ACFO^-$-model associated to a model  $(F, \triangleleft)$ of $\ACFO^-$, and suppose  $(F^\times, \Sigma, \triangleleft)$ is existentially closed. We assume further that $(F, \triangleleft) \models \ACFO_p^-$  for a fixed $p$ and omit the proof of the similar case where $(F, \triangleleft)$ is a model of $\ACFO_0^-$. Note that an $\ACFO^-$-model extending $(F,  \triangleleft)$ is then automatically an  $\ACFO_p^-$-model. 
Let $\phi(x)$ be a quantifier-free $\lac$-formula which defines a nonempty set in some $\ACFO_p^-$-model extending $(F,  \triangleleft)$. To get the backward implication, we need to show that $\phi(x)$  already defines a nonempty set in $(F,  \triangleleft)$. 

Note that $\phi(x)$  is logically equivalent to $(  \phi(x) \wedge x_i =0) \vee (\phi(x) \wedge x_i \neq 0)$ with $i \in \{1, \ldots, m\}$, and $\phi(x) \wedge x_i =0$ is equivalent over $\ACFO_p^-$ to a quantifier-free formular with fewer variables. So we reduced to the case where $\phi(x)$ is logically equivalent to  $\phi(x) \wedge \bigwedge_{i=1}^m x_i \neq 0$. Consider the special case where the only atomic formulas of $\phi(x)$ in which $+$ appears are of the form 
$$t_1(x) + \ldots +t_n(x) =t_{n+1}(x)$$
where  $t_i(x)$ does not further contain $+$ for $i \in \{1, \ldots, n+1\}$. Such a formula $t_1(x) + \ldots +t_n(x) =t_{n+1}(x)$ defines in an arbitary $\ACFO_p^-$-model  the same set that $\Sigma_{n+1}(t_1(x) , \ldots ,t_n(x), t_{n+1}(x))$ defines in the associated punctured $\ACFO_p^-$-model. So we get an $\lac^\times$ formula $\phi^\times(x)$ such that $\phi(x)$ defines in an $\ACFO_p^-$-model  the same set that $\phi^\times(x)$ defines in its associated punctured $\ACFO_p^-$-model. In particular,  $\phi^\times(x)$ defines a nonempty set in a punctured $\ACFO_p^-$-model extending $(F^\times, \Sigma, \triangleleft)$. As $(F^\times, \Sigma, \triangleleft)$ is existentially closed, $\phi^\times(x)$ defines a nonempty set in $(F^\times, \Sigma, \triangleleft)$.  Thus, in this special case,  $\phi(x)$ defines a nonempty set in $(F,  \triangleleft)$ .

We next reduce to the special case in the preceding paragraph. Arrange that every term $t(x)$ appearing in $\phi(x)$ is of the form $t_1(x) + \ldots +t_k(x)$ where  $t_i(x)$ does not contain $+$ for $i \in \{1, \ldots, k\}$.  Introduce a new variable $y_t$ for each such $t$, and set $y$ to be the tuple of variables built up from such $y_t$. Replace the appearance of each aforementioned $t(x)$ in $\phi(x)$ with $y_t$ to get a formula $\psi(x, y)$. Then $\phi(x)$ is equivalent across $\ACFO_p^-$-models to $$  \exists y \big(  \psi(x,y) \wedge  \bigwedge_t t(x) =y_t       \big).  $$
Note that $ \psi(x,y) \wedge  \bigwedge_t y_t = t(x)$ is of the form in the preceding paragraph and defines a nonempty set in an $\ACFO_p$-models extending $(F,  \triangleleft)$. So $ \psi(x,y) \wedge  \bigwedge_t y_t = t(x)$ defines a nonempty set in $(F,  \triangleleft)$. Hence,  $\phi(x)$ also defines a nonempty set in $(F,  \triangleleft)$, which concludes the proof of the backward implication.

The forward implication is much easier. The main difficulty with the backward implication comes from the fact that $+$ has no corresponding function symbol in $\lac^\times$. On the other hand, basic functions of  $(F^\times, \Sigma, \triangleleft)$  are restrictions of basic functions of $(F, \triangleleft)$ to $0$-definable sets of $(F, \triangleleft)$, and basic relations of $(F^\times, \Sigma, \triangleleft)$ are $0$-definable sets of $(F, \triangleleft)$. So the analogous argument can be carried out without worrying about the aforementioned difficulty.
\end{proof}

\noindent In order to ``geometrically'' characterize the existentially closed models of the expansion of a theory by a unary predicate, Chatzidakis and Pillay implicitly introduced a ``largeness property'' for definable sets \cite{ChatPillay}; the name largeness here is taken from \cite{ShelahSimon}. We will provide an analogous notion in  this setting.

\meno Suppose $(G, \square)$ is an expansion of $G$,  $X \subseteq G^m$ is definable in $(G, \square)$, and $(\mathbf{G}, \square)$ is a monster elementary extension of $(G, \square)$. We say that $X$ is { \bf multiplicatively large} if there is $a \in X(\mathbf{G})$ which is not a solution of any nontrivial multiplicative equations over $G$. 

\meno The above definition in particular applies to definable sets in a circularly ordered abelian group
$(G, \triangleleft)$ and definable sets in a punctured field $(G, \Sigma)$. We can extend it in an obvious way to cover definable sets in a field $F$. A definable subset $X \subseteq F^m$ is {\bf multiplicative large} if $X \cap (F^\times)^m$ is multiplicatively large as a definable subset of the punctured field $(F^\times, \Sigma)$ associated to $F$.
 
\begin{rem} \label{Rem: Remarks on multiplicative largeness}
Suppose $(G, \square)$ is an expansion of $G$, and $X$, $X_1$, and $X_2$ are definable in $(G, \square)$ with $X = X_1 \cup X_2$. Then $X$ is multiplicatively large if and only if  either $X_1$ is multiplicatively large or $X_2$ is multiplicatively large.
\end{rem}

\noindent Even though  multiplicative largeness can be defined very generally, it only behaves well under stronger assumptions:

\begin{lem} \label{lem: Preservation of multiplicative largeness under elementary extension}
Suppose $(G, \square)$ is an expansion of a model of $\thgg$, and $X \subseteq G^m$ is a multiplicatively large definable set in $(G, \square)$. 
 If  $(G', \square)$ is an elementary extension of $(G, \square)$, then $X(G')$ is a multiplicatively large definable set in $(G', \square)$.
\end{lem}
\begin{proof}
Suppose $(G, \square)$, $(G', \square)$,  and $X$ are as above. It follows from Lemma~\ref{dimensionandmultiplicativelarge} that $\text{mdim}(X) = m$. As $\text{mdim}$ coincides with Morley rank which is preserved under taking elementary extension, $\text{mdim}(X(G')) =m$.  Applying  Lemma~\ref{dimensionandmultiplicativelarge} again, we get that $X(G')$ is also multiplicatively large.
\end{proof}

\noindent We now move on to give the geometric characterization we promised.  Suppose $(G, \Sigma, \triangleleft)$ is a punctured $\ACFO_p^-$-model. We say that $(G, \Sigma, \triangleleft)$ satisfies the {\bf geometric characterization} if the following two conditions are satisfied:
\begin{enumerate}
    \item $(G, \triangleleft) \models \thmm$;
    \item if $X_1 \subseteq G^m$ is a multiplicatively large qf-set in $(G, \Sigma)$ and $X_2 \subseteq G^m$ is a multiplicative large qf-set in $(G, \triangleleft)$, then $X_1 \cap X_2 \neq \emptyset$.
\end{enumerate}
We say that an $\ACFO^-$-model $(F, \triangleleft)$ satisfies the {\bf geometric characterization} if its associated punctured $\ACFO^-$-model does, or equivalently, if it satisfies a definition as above but with $(G, \triangleleft)$ replaced by $(F^\times, \triangleleft)$ and $(G, \Sigma)$ replaced by $F$. Note that $\ACF^\times$ admits quantifier elimination, so the assumption that $X_1$ is a qf-set is in fact unnecessary.

\meno In the rest of the section, we will show that an $\ACFO^-$-model is existentially closed if and only if it satisfies the geometric characterization. If $T$ is an $L$-theory, let $T(\forall)$ denote the set of $L$-consequences of $T$.

\begin{lem} \label{Lem: Enriched embedding 1}
Suppose $(G, \Sigma, \triangleleft)$ is an $\lac^\times$-structure with $(G, \Sigma) \models \ACF_p^\times (\forall)$ and $(G, \triangleleft) \models \thmp(\forall)$. Then  $(G, \Sigma, \triangleleft)$ can be $\lac^\times$-embedded into a punctured $\ACFO_p^-$-model $(G', \Sigma, \triangleleft)$ such that $(G', \triangleleft)$ is a model of $\thmm$. A similar statement holds for an $\lac^\times$-structure $(G, \Sigma, \triangleleft)$ with $(G, \Sigma) \models \ACF_0^\times (\forall)$ and $(G, \Sigma) \models \thmo(\forall)$.
\end{lem}

\begin{proof}
We will only prove the first statement as the proof of the second statement is similar. Let $(G, \Sigma, \triangleleft)$ be as above.
We will construct a sequence $(G_n, \Sigma, \triangleleft)_n$ of $\lac^\times$-structures such that 
\begin{enumerate}
    \item $(G_0, \Sigma, \triangleleft)$ extends  $(G, \Sigma, \triangleleft)$ as an $\lac^\times$-structure;
    \item $(G_{n+1}, \Sigma, \triangleleft)$ is an $\lac^\times$-extension of $ (G_n, \Sigma, \triangleleft)$;
    \item $(G_{2n}, \Sigma)$ is a punctured $\ACFO_p^-$-model;
    \item $(G_{2n+1}, \triangleleft) \models \thmp$ and $(G_{2n+1}, \Sigma)$ is a model of $\ACF_p^\times(\forall)$.
\end{enumerate}
Then we can take $(G', \Sigma, \triangleleft)$ to be the union of $(G_n, \Sigma, \triangleleft)_n$. Note that $\ACF^\times_p$ and $\thmp$ are both inductive theories, as they are model complete, so it is easy to see that $(G', \Sigma, \triangleleft)$ satisfied the desired conclusion.

  As $(G, \Sigma) \models \ACF_p^\times (\forall)$, we can get $(G_0, \Sigma)\models \ACF_p^\times$ extending $(G, \Sigma)$ as an $\lart$-structure. On the other hand $(G, \triangleleft) \models \thmp(\forall)$, so we can get a monster model $(\mathbf{G}, \triangleleft)$ of $\thmp$ extending $(G, \triangleleft)$ as an $\lam$-structure. Now both $G_0$ and $\mathbf{G}$ are models of $\thgp$, and this theory admits quantifier elimination by Proposition~\ref{prop: GMproperties}. Hence, there is an embedding of $G_0$ into $\mathbf{G}$. We can then define $\triangleleft$ on $G_0$ by pull-back via $f$.  Clearly, $(G_0, \Sigma, \triangleleft)$ is a model of $\ACFO^-_p$.
  
  Suppose we have constructed  $(G_{2n}, \Sigma, \triangleleft)$ satisfying both (2) and (3).  Then $(G_{2n}, \triangleleft) $ is a model of $\thmp^-$ which has $\thmp$ as a model companion, so we can get $(G_{2n+1}, \triangleleft) \models \thmp$ extending $(G_{2n}, \triangleleft)$ as an $\lam$-structure. Choose a monster model $(\mathbf{G}, \Sigma)$ of $\ACF_p^\times$ extending $(G_{2n}, \Sigma)$. Note that $G_{2n}$, $G_{2n+1}$, and $\mathbf{G}$ are models of $\thgp$, which is a model complete $\lag$-theory. So there is an $\lag$-embedding $f: G_{2n+1} \to \mathbf{G}$ which extends the identity map on $G_{2n}$. Define the family of relations $\Sigma$  on $(G_{2n+1}, \triangleleft)$ as the pull-back via $f$ of the family $\Sigma$ on $\mathbf{G}$. Then $G_{2n+1} \models \ACF_p^\times(\forall)$ by construction, and so $(G_{2n+1}, \triangleleft) \models \thmp$  satisfies (4).

Finally, suppose we have constructed  $(G_{2n+1}, \Sigma, \triangleleft)$ satisfying both (2) and (4).  We note that the only thing used in the preceding paragraph is that $\thmp$ is the model companion of $\thmp^-$ and that $\thgp$ is model complete. Hence, we can carry out exactly the same strategy to get the desired conclusion by replacing the former with the fact that $\ACF^\times_p$ is the model companion of $\ACF_p(\forall)$ and reusing the latter. 
\end{proof}

\begin{rem}
Using the fact that $\thgp$ is strongly minimal, one can produce a quicker proof of Lemma~\ref{Lem: Enriched embedding 1} by constructing $(G', \Sigma, \triangleleft)$ directly. We still choose to present the longer proof here to make the neccesary ingredients transparent.
\end{rem}

\noindent We need another embedding lemma:

\begin{lem} \label{lem: enriched embedding 2}
Suppose $(G, \Sigma, \triangleleft)$ is a punctured $\ACFO_p^-$-model with $(G, \triangleleft) \models \thmp$. Then  $(G, \Sigma, \triangleleft)$ can be $\lac^\times$-embedded into a punctured $\ACFO_p^-$-model which satisfies the geometric characterization. A similar statement holds for a punctured $\ACFO_0^-$-model $(G, \Sigma, \triangleleft)$ with $(G, \triangleleft) \models \thmo$.
\end{lem}

\begin{proof}
We will only prove the first statement, as the proof for the second statement is similar. Suppose $(G, \Sigma, \triangleleft)$ is as stated. Let $ X_1 \subseteq (F^\times)^m$ is a multiplicative large qf-set in $(G, \Sigma)$ and  $ X_2 \subseteq (F^\times)^m$ is a multiplicative large qf-set in $(G, \triangleleft)$. Our problem can be reduced to finding $(G', \Sigma, \triangleleft) \models \ACFO_p^-$ extending  $(G, \Sigma, \triangleleft)$ with $(G', \triangleleft) \models \thmm$ such that $X_1(G') \cap X_2(G') \neq \emptyset$. Indeed, we can simply iterate this construction and take the union.

We now construct the aforementioned $(G', \Sigma, \triangleleft)$. Take a monster elementary extension $(\mathbf{G_1}, \Sigma)$ of $(G, \Sigma)$ and a monster elementary extension $(\mathbf{G_2}, \triangleleft)$ of $(G, \triangleleft)$. As $X_1$ is multiplicatively large, we get $a' \in X_1(\mathbf{G}_1)$ whose components are multiplicatively independent over $G$. Likewise, we get $b' \in X_1(\mathbf{G}_2)$ whose components are multiplicatively independent over $G$. 
Let $f:  \la G, a' \ra \to \mathbf{G}_2$ be the unique map which extends the identity map on $G$ and maps $a'$ to $b'$. Define the family $\Sigma$ on $\la G, a' \ra$ by restricting the family with the same name on $\mathbf{G_1}$, and define the relation $\triangleleft$ on $\la G, a' \ra$ by pulling back via $f$ the relation with the same name on $\mathbf{G_2}$. Then $( \la G, a'\ra, \Sigma, \triangleleft) $ is an $\lac^\times$-structure with $( \la G, a'\ra, \Sigma) \models \ACF^\times(\forall)$ and $( \la G, a'\ra, \Sigma) \models \thmm(\forall)$. Applying Lemma~\ref{Lem: Enriched embedding 1} yields the desired conclusion.
\end{proof}

\noindent  Lemma~\ref{lem: enriched embedding 2} essentially gives us that the characterization works in one direction.

\begin{cor} \label{existentiallyclosedimpliesgeneric}
Suppose $(G, \Sigma, \triangleleft)$ is an existentially closed punctured $\ACFO^-$-model. Then $(G, \triangleleft)$ satisfies the geometric characterization. 
\end{cor}
\begin{proof}
We consider only the case where $(G, \Sigma, \triangleleft)$ is a punctured $\ACFO_p^-$-model; the other case with  $(G, \Sigma, \triangleleft)$ a punctured $\ACFO_0^-$-model is very similar.
Using Lemma~\ref{lem: enriched embedding 2}, we obtain a punctured $\ACFO_p^-$-model  $(G', \Sigma, \triangleleft) $ extending $(G, \Sigma, \triangleleft)$ as an $\lac^\times$-structure such that $(G', \Sigma, \triangleleft) $  satisfies the geometric characterization. 

The structure $(G, \Sigma, \triangleleft)$ is existentially closed in $(G', \Sigma, \triangleleft)$, so $(G, \triangleleft)$ is existentially closed in $(G', \triangleleft)$. As $(G', \Sigma, \triangleleft) $ satisfies the  geometric characterization, $(G', \triangleleft)$ is a model of $\thmm_p$. The theory $\thmm_p$ is model complete, so we can assume that it only consists of $\forall \exists$ statements. It follows that $(G, \triangleleft)$ is also a model of $\thmm_p$.

Suppose $X_1 \subseteq G^m$ is a multiplicatively large qf-set in $(G,  \Sigma)$ and $X_2 \subseteq G^m$ is multiplicatively large qf-set in $(G, \triangleleft)$. 
Then $X_1(G' )$ is multiplicatively large in $(G',  \Sigma)$, and  $X_2(G')$ is multiplicatively large in  $(G', \triangleleft)$ by Lemma~\ref{lem: Preservation of multiplicative largeness under elementary extension} and the fact that both $\ACF_p^\times$ and $\thmp$ are model complete. As $(G', \Sigma, \triangleleft) $  satisfies the geometric characterization, $X_1(G') \cap  X_2(G') \neq \emptyset$. So $X_1 \cap X_2 \neq \emptyset$ as well by the fact that $(G, \Sigma, \triangleleft)$ is existentially closed. The desired conclusion follows.
\end{proof}

\noindent We next verify that the characterization works in the other direction:

\begin{lem} \label{genericimpliesexistentiallyclosed}
Suppose $(G, \Sigma, \triangleleft)$ is a punctured model of $\ACFO^-$ that satisfies the geometric characterization. Then $(G, \Sigma, \triangleleft)$ is existentially closed.
\end{lem}

\begin{proof}
It suffices to prove the corresponding statements for punctured $\ACFO_p^-$-models and punctured  $\ACFO_0^-$-models. We will only prove the former as the latter is very similar.
 Suppose $(G, \Sigma, \triangleleft)$ is a generic punctured model of $\ACFO_p^-$, and $\phi(x)$ is a quantifier-free $\lac^\times(G)$-formula which defines a nonempty set in a punctured $\ACFO_p^-$-model $(G', \Sigma, \triangleleft)$ extending $(G, \Sigma, \triangleleft)$. Our job is to show that $\phi(x)$ defines  a nonempty set in $(G, \Sigma, \triangleleft)$. 

As the only function symbols in both $\lart$ and $\lam$  already appear in $\lag$, we can reduce to the case where $\phi(x) = \psi(x) \wedge \theta(x)$ with $\psi(x)$ a quantifier-free $\lart(G)$-formula and $\theta(x)$ a quantifier-free $\lam(G)$-formula. Let $a' \in (G')^m$ be such that $(G', \Sigma, \triangleleft) \models \phi(a')$. Suppose $a'$ is multiplicatively independent. Then  $\psi(x)$ and $\theta(x)$ define multiplicatively large sets in $(G, \Sigma)$ and $(G, \triangleleft)$.  So  $\phi(x)$ defines a nonempty set in $(G, \Sigma, \triangleleft)$ by the assumption that $(G, \Sigma, \triangleleft)$ satisfies the geometric characterization. 

Now consider the general case where $a'$ might not be multiplicatively independent. Then we can choose a tuple $b' \in (G')^n$ which is multiplicatively independent such that $a'= t(b')$ with $t =(t_1, \ldots, t_m)$ and $t_i(y)$ is an $\lag$-term for $i \in \{1, \ldots, m\}$. Applying the previous case for the formula  $\phi(t(y)) = \psi(t(y)) \wedge \theta(t(y))$, we obtain $b \in G^y$ such that $(G, \Sigma, \triangleleft) \models \phi(t(b))$. Thus,  $\phi(x)$ defines in  $(G, \Sigma, \triangleleft)$ a nonempty set, which is our desired conclusion.
\end{proof}

\noindent Finally, we put everything together:

\begin{prop} \label{prop: generic and existentially closed}
An $\ACFO^-$-model $(F, \triangleleft)$ is existentially closed if and only if $(F, \triangleleft)$ satisfies the geometric characterization.
\end{prop}
\begin{proof}
This follows easily from Lemma~\ref{Existentially closed models exchange}, Corollary~\ref{existentiallyclosedimpliesgeneric}, and  Lemma~\ref{genericimpliesexistentiallyclosed}.
\end{proof}

\subsection{Axiomatization} \label{Sec: Axiomatization}
Just as in \cite{ChatPillay}, we want to establish that the class of models of $\ACFO^-$ satisfying the geometric characterization is elementary. In order to do so, 
the key is to show that $\ACF$ and $\thmm$ each ``defines multiplicative largeness''.
 
\meno Throughout Section~\ref{Sec: Axiomatization},  $F$ is an  algebraically closed field, and $V \subseteq (F^\times)^m$ is a quasi-affine variety (i.e, a Zariski-open subset of an irreducible Zariski-closed subset of $(F^\times)^m$). We equip $(F^\times)^m$ with the group structure  given by coordinate-wise multiplication, and let $1^{(m)}$ be the identity element of $(F^\times)^m$.

\meno Suppose $T$ is an $L$-theory, $\phi(x,y)$ is an $L$-formula, and $\sP$ is a property of $\sM$-definable subsets of $M^x$ with $\sM \models T$. We say that $T$ {\bf defines $\sP$ for $\phi(x,y)$} if there is an  $L$-formula $\delta(y)$ such that for all $\sM \models T$ and $b \in M^y$, we have
$$ \phi(\sM, b) \text{ satisfies } \sP \quad \text{ if and only if} \quad \sM \models \delta(b).         $$
A formula  $\delta(y)$ as above is called a {\bf $\sP$-defining formula} over $T$ for $\phi(x,y)$. The fact that this is key to axiomatization can be seen through the following lemma:
\begin{lem} \label{Multiplicative largeness and axiomatization}
Assume $\ACF$ defines multiplicative largeness for all $\lar$-formulas $\phi(x, y)$ and $\thmm$ defines multiplicative largeness for all quantifier-free $\lam$-formulas $\psi(x, z)$. Then the class of $\ACFO^-$-models satisfying the geometric characterization is elementary.
\end{lem}
\begin{proof}
Let $\ACFO^{-1 / 2}$ be the $\lac$-theory whose models  are the $(F, \triangleleft) \models \ACFO^-$ with $(F^\times, \triangleleft) \models \thmm$. Obtain $\ACFO$ from $\ACFO^{-1 / 2}$ by  adding for each $\lar$-formula $\phi(x,y)$ and $\lam$-formula $\psi(x,z)$ the formula
$$ \forall y \forall z\big(  \delta(y) \wedge \hat{\theta}(z) \to \exists x ( \phi(x, y) \wedge \hat{\psi}(x,z)) \big), $$
where $\delta(y)$ is a multiplicative largeness-defining formula over $\ACF$ for $\phi(x,y)$,  $\theta(z)$ is a multiplicative largeness-defining formula over $\thmm$ for $\psi(x,z)$, and $\hat{\theta}(z)$ and $\hat{\psi}(x,z)$ are the obvious modifications of $\theta(z)$ and $\psi(x,z)$ such that $\hat{\theta}(z)$ and $\hat{\psi}(x,z)$  apply to all tuples with components in $F$. It is easy to see that $\ACFO$ axiomatizes the class  of $\ACFO^-$-models satisfying the geometric characterization.
\end{proof}

\noindent We next obtain various characterizations of multiplicative largeness in models of $\ACF$ and deduce from them the first condition of Lemma~\ref{Multiplicative largeness and axiomatization}. 

\begin{lem} \label{Multiplicatively large for variety}
A quasi-affine variety $V \subseteq (F^\times)^m$ is multiplicatively large if and only if no nontrivial multiplicative equation vanishes on $V$.
\end{lem}
\begin{proof}
The forward implication is immediate. Suppose  no nontrivial multiplicative equation vanishes on $V$. As $V$ is irreducible, it is not a subset of a finite union of solution sets of nontrivial multiplicative equations. The desired conclusion then follows a compactness argument. 
\end{proof}

\noindent  If $B \subseteq F^\times$, a
{\bf multiplicative system} over $B$ is simply a conjunction of multiplicative equations over $B$.
Fact~\ref{MultLarge0} below about definable subgroups of $(F^\times)^m$ is a consequence of the fact that definable subgroups of $(F^\times)^m$ are closed and of the characterization of algebraic subgroups of $(F^\times)^m$. For the former, see for instance~\cite[Lemma 7.4.9] {Marker}. For the latter, see for instance~\cite[Corollary 3.2.15]{Bombieri}; the proof there is for characteristic $0$ but goes through in positive characteristics.  

\begin{fact} \label{MultLarge0}
Every connected definable subgroup of $(F^\times)^m$ is defined by a multiplicative system over $\emptyset$.
\end{fact}

\noindent Suppose $X_1, \ldots, X_n$ are subsets of $G$. Set $X_1\cdots X_n = \{a_1 \cdots a_n \mid a_i \in X_i  \text{ for }   1\leq i \leq n\}.$ Moreover, if $X_1 = \dots = X_n =X$, then we denote this as  $\Pi_n(X) $.  The following fact is a special case of Zilber's Indecomposability theorem but was  also known much earlier; see Theorem 7.3.2 of \cite{Marker}.
\begin{fact} \label{Zilbertheorem}
Suppose $1^{(m)}$ is in $V$.  Then $\Pi_{2m}(V)$ is a connected definable subgroup of $(F^\times)^m$. Hence, $\Pi_{2m}(V)$ is a subgroup of every definable subgroup of $(F^\times)^m$  containing $V$ as a subset.
\end{fact}

\noindent We now have a simple criterion for multiplicative largeness:

\begin{lem} \label{MultLarge1}
If $1^{(m)}$ is in $V$, then $V$ is multiplicatively large if and only if $\Pi_{2m}(V)$ is $(F^\times)^m$.
\end{lem}

\begin{proof}
For the forward implication, suppose $V$ is multiplicatively large. Then by Lemma~\ref{Multiplicatively large for variety} and Fact~\ref{MultLarge0}, no proper definable subgroup of $(F^\times)^m$ contains $V$ as a subset. It then follows from Fact~\ref{Zilbertheorem} that $\Pi_{2m}(V) =  (F^\times)^m$. The backward implication is immediate from Lemma~\ref{Multiplicatively large for variety}.
\end{proof}

\noindent To get from quasi-affine varieties over $F$ to general sets definable in $F$, we need a result related to defining irreducibility. This and other related results are included in Fact~\ref{fact: Definable of irreducibility} as we will also need them later on; see \cite[Chapter 10]{Will} for details. 

\begin{fact} \label{fact: Definable of irreducibility}
Supose $\phi(x, y)$ is an $\lar$-formula, $d$ is in $\nn$, and $r$ is in $\nn^{\geq 1}$. Then there are formulas $\delta_d(y)$, $\mu_r(y)$,  $\iota(y)$, and $\psi(x,z)$ such that if the families $(X_b)_{b\in Y}$ and $(X_c)_{c\in Z}$ are defined  in $F$
by $\phi(x, y)$ and $\psi(x,z)$, then we have the following:
\begin{enumerate}
    \item[\ir] $F \models \delta_d(b)$ for $b\in Y$ if and only if $\dim(X_b) =d$;
    \item[\iir] $F \models \mu_r(b)$ for $b\in Y$ if and only if the Morley degree of $X_b$ is $r$;
    \item[\iir] $F \models \iota(b)$ for $b\in Y$ if and only if $X_b$ is a quasi-affine subvariety of $F^m$;
    \item[\ivr] $(X_c)_{c\in Z}$ is a family of quasi-affine varieties which contains all irreducible components of members of $(X_b)_{b\in Y}$.
\end{enumerate}
\end{fact}

\noindent We now put together Fact~\ref{fact: Definable of irreducibility} and Lemma~\ref{MultLarge1}:

\begin{prop} \label{MultLarge2}
The theory $\ACF$ defines multiplicative largeness for every $\lar$-formula $\phi(x,y)$.
\end{prop}
\begin{proof}
Suppose $\phi(x,y)$ is an arbitrary $\lart$-formula, $\psi(x, z)$ is as in Fact~\ref{fact: Definable of irreducibility}, and $\phi(x,y)$ and $\psi(x, z)$ define in $F$ the families $(X_b)_{b \in Y}$ and $(X_c)_{c \in Z}$. Observe that $X_b$ with $b \in Y$ is multiplicatively large if and only if there is $c \in Z$ with $X_c \subseteq X_b$ and $a \in X_c \cap (F^\times)^m$ such that
$\Pi_{2m}\big(a^{-1}(X_c \cap (F^\times)^m)\big)=(F^\times)^m$. It is easy to see from here that $\ACF$ defines multiplicative largeness for $\phi(x,y)$.
\end{proof}

\noindent Every multiplicative equation over $\emptyset$ is equivalent over $\thgg$ to a multiplicative equation $t(x) =t'(x)$ of the simplified form  where the power of every variable is nonnegative and each variable appears only on at most one side of the equation. The {\bf degree} of a multiplicative equation is the highest power of $x_i$ which appears in the above simplified form as $i$ ranges over $\{1, \ldots, m\}$. A simple application of the compactness theorem gives the following corollary, which we will use later on.

\begin{cor} \label{MultLarge3}
Suppose  $(V_b)_{b\in Y}$ is a family of quasi-affine subvarieties of $(F^\times)^m$ passing through $1^{(m)}$. Then there is $N >0$ such that for all $b\in Y$, either $V_b$ is multiplicatively large or a nontrivial multiplicative equation over $\emptyset$ with degree at most $N$ vanishes on $V_b$.
\end{cor}

\noindent We will next obtain a characterization of multiplicative largeness in models of $\thmm$ and from there obtain the second condition of Lemma~\ref{Multiplicative largeness and axiomatization}.

\meno Suppose $(G, \triangleleft) \models \thmm$. The {\bf $\triangleleft$-topology} on $G^m$ is defined as the topology which has a basis consisting of sets of the form $U = U_1 \times \cdots \times U_m$ where $U_i$ is the ``interval''
$$ \{a \in G : \triangleleft(d_i, a, d'_i) \}$$
with $d_i$ and $d'_i$ in $G$ for $i \in \{1, \ldots, m\}$. It is also easy to see that the  $\triangleleft$-topology on $G^m$ is simply the product of the $\triangleleft$-topologies on the $m$ copies of $G$. 
\begin{lem} \label{opensetandcontinuosmap}
Suppose $(G, \triangleleft)$ is a model of $\thmm$. Then we have the following:
\begin{enumerate}
    \item [\ir] $\triangleleft$ as a subset of $G^3$ is $\triangleleft$-open;
    \item [\iir] the multiplication map is continuous.
\end{enumerate}
\end{lem}
\begin{proof}
It is well known that $\uu_{(p)}$ for an arbitrary $p$ and $\uu$ are dense in $\mathbb{T}$ with respect to the Euclidean topology. Hence, when $(G, \triangleleft)$ is $(\uu_{(p)}, \triangleleft)$ for some given $p$ or $(\uu, \triangleleft)$, the $\triangleleft$-topology is just the subspace topology with respect to the usual Euclidean topology on $\mathbb{T}$. Hence, (i) and (ii) are automatic in these cases. Note that properties (i) and (ii) can be expressed as $\lam$-statements, so the desired conclusion follows from Proposition~\ref{Up is model of GMO-} and  Proposition~\ref{prop: GMOproperties}.
\end{proof}

\begin{prop} \label{prop: multiplicative largeness in GMO}
Suppose $(G, \triangleleft)$ is a model of $\thmm$ and $X \subseteq G^m$  is defined in $(G, \triangleleft)$ by a quantifer-free $\lam(G)$-formula. Then $X$ is multiplicatively large if and only if  $X$ contains a nonempty subset which is  $\triangleleft$-open in $G^m$.
\end{prop}
\begin{proof}
Let $(G, \triangleleft)$ and $X$ be as stated above. For the forward implication, suppose $X$ is multiplicative large. Note that $\neg \triangleleft(x,y,z)$ is equivalent over $\thmm$ to $$\triangleleft(z, y, x) \vee (x=y) \vee(y=z) \vee(z=x).$$ So quantifier free $\lam$-formulas are equivalent over $\thmm$ to positive quantifier-free formulas. Applying also Remark~\ref{Rem: Remarks on multiplicative largeness}, we can assume that $X$ is defined by a formula of the form
$$ \bigwedge_{i \in I} \triangleleft\big( t_i(x), t'_i(x), t''_i(x)\big) \wedge \bigwedge_{j \in J} \big(t_j(x) =t'_j(x)\big). $$
where $t_i(x), t'_i(x),t''_i(x), t_j(x), t'_j(x)$ are $\lag(G)$-terms for all $i \in I$ and $j \in J$. 
As $X$ is multiplicatively large, one must have that $t_j(x) = t'_j(x)$ is trivial for all $j \in J$.  It then follows from Lemma~\ref{opensetandcontinuosmap} that $X$ is open, which gives us the desired conclusion.

For the backward implication, suppose $X$ contains a a nonempty subset which is  $\triangleleft$-open in $G^m$. We may assume that $X = \prod_{i =1}^m X_i$ where 
$$X_i
= \{ a_i \in G \mid  \triangleleft(d_i, a_i, d'_i) \}, $$
with $d_i, d'_i \in G$ and $d_i \neq d'_i$ for $i \in \{1, \ldots, m\}$. Let $(\mathbf{G}, \triangleleft)$ be a monster elementary extension of $(G, \triangleleft)$. Using Proposition~\ref{prop: GMOproperties}, it is easy to show that $X_i$ is infinite by reducing to the special cases where $G= \uu$ or $G= \uu_{(p)}$ for some $p$. Hence, $|X_i(\mathbf{G})| > |G|$ for $i \in I$. Hence, we can choose the desired $a' = (a'_1, \ldots, a'_m)$ in $ X(\mathbf{G})$ by ensuring that $a'_{i+1} \in X_i( \mathbf{G}) $ is multiplicatively independent over $a'_1, \ldots, a'_{i}$ for $i \in \{1, \ldots, m-1\}$.
\end{proof}
\noindent From the definition of $\triangleleft$-topology, we immediately get:
\begin{cor} \label{GMO defines multiplicative largeness}
The theory $\thmm$ defines multiplicative largeness for all quantifier free $\lam$-formulas $\phi(x,y)$.
\end{cor}
\noindent Combining Proposition~\ref{prop: multiplicative largeness in GMO},  Remark~\ref{Rem: Remarks on multiplicative largeness}, and the well-known fact that definable sets in $\ACF$-models are finite unions of quasi-affine varieties, we get:
\begin{cor} \label{refineddefinition of generic}
Suppose $(F, \triangleleft)$ is a model of $\ACFO^-$ with  $(F^\times, \triangleleft) \models \thmm$. Then $(F, \triangleleft)$ satisfies the geometric characterization if and only if 
all multiplicatively large $V\subseteq (F^\times)^m$ are dense with respect to the $\triangleleft$-topology.
\end{cor}

\noindent We now put together the results of this section and the preceding section to get the existence of the model companion $\ACFO$ of $\ACFO^-$.

\begin{proof}[Proof of Theorem~\ref{thm: Model companion of ACFO-}] The desired conclusion follows immediately from Proposition~\ref{prop: generic and existentially closed}, Lemma~\ref{Multiplicative largeness and axiomatization}, Proposition~\ref{MultLarge2}, and Corollary~\ref{GMO defines multiplicative largeness}.
\end{proof}

\noindent As a side note, we will show that every model of $\ACFO$ has 
 $\text{TP}_2$. The notion was defined in \cite{ShelahdefinitionofNTP2} by Shelah and systematically studied in \cite{Chernikov} by Chernikov. 
We use here a finitary version of the definition given in  \cite{Chernikov}. Let $\sM$ be an $L$-structure. An  $L$-formula $\phi(x,y)$ witnesses that $\sM$ has $\mathrm{TP}_2$ if for each finite set $I$, there is a family  $(b_{ij})_{(i,j) \in I^2}$ of elements of $M^n$ such that the following conditions hold:
\begin{enumerate}
    \item $\sM \models \neg \exists x \big(\phi(x, b_{ij}) \wedge \phi(x,b_{ij'}) \big)$ for every $i \in I$ and distinct $j$ and $j'$ in $I$;
    \item $\sM \models \exists x \bigwedge_{i \in I}\phi(x, b_{if(i)})$  for any $f: I \to I$.
\end{enumerate}
We say that $\sM$ has $\mathrm{TP}_2$ if there is a formula $\phi(x,y)$ which witnesses that $\sM$ has $\mathrm{TP}_2$ and say that $\sM$ has $\mathrm{NTP}_2$ otherwise. 

\begin{prop} \label{Prop: this structure has TP_2}
Every model of $\ACFO$ has $\mathrm{TP}_2$. 
\end{prop}

\begin{proof}
Suppose  $(F, \triangleleft)$  is a model of $\ACFO$. Let $x$ be a single variable, $y = (z, t, t')$ with $z$, $t$, and $t'$ single variables, and $\phi(x,y)$ the formula
$$  \triangleleft(x+ z, t, t').  $$
Let $I$ be an arbitrary finite set. Get a family $(c_i)_{i \in I}$ of distinct elements in $F$. Obtain a family  $(d_j, d'_j)_{j \in I}$ of pairs of elements in $F^\times$ with $d_j \neq d'_j$ for all $j \in J$ and   
 $$ (F, \triangleleft) \models \neg\exists x \big(\triangleleft(d_j, x, d'_j) \wedge \triangleleft(d_{j'}, x, d'_{j'})\big) \quad \text{ for distinct } j, j' \in I. $$ Set $b_{ij} = (c_i, d_j, d'_j)$.
 It is easy to see that $\phi(x,y)$ together with $(b_{ij})_{(i,j) \in I^2}$ satisfy (1) in the definition of a $\mathrm{TP}_2$-witness. We assume without loss of generality that $I = \{1, \ldots, k\}$. Set
$$V = \{(a+c_1, \ldots, a+c_k) : a \in F\}.$$ It suffices to show that $V$ is multiplicatively large as it will follow that $\phi(x,y)$ together with $(b_{ij})_{(i,j)\in I^2}$ satisfies (2) in the definition of a $\mathrm{TP}_2$-formula as well.
 We can reduce further to showing triviality for an arbitrarily chosen multiplicative equation $x_1^{n_1}\cdots x_m^{n_m} = c x_1^{n'_1}\cdots x_k^{n'_m}$ vanishing on $V$ where $c$ is in $F^\times$, and  $n_i$ and $n'_i$ are in $\nn$ with either $n_i=0$ or $n'_i=0$ for $i \in \{1, \ldots, m\}$. In this case, we have
$$ (a+c_1)^{n_1}\cdots (a+c_m)^{n_m} = c (a+c_1)^{n'_1}\cdots (a+c_m)^{m'_m} \quad \text{for all } a \in F. $$
  For $i \in \{1, \ldots, m\}$, we substitute $a = -c_i$  and deduce $n_i =n'_{i}=0$. Hence, we also get $c=1$, and the desired conclusion follows.
\end{proof}

\begin{rem}
Proposition~\ref{Prop: this structure has TP_2} is surprising: following  \cite{ChatPillay}, one would expect that models of $\ACFO$ have $\NTP_2$. It suggests that $\text{NTP}_2$ is not quite 
``stable $+$ order $+$ random''.
In the same direction, recent evidence seems to suggest that ``stable $+$ random'' is $\text{NSOP}_1$ instead of simple as earlier thought~\cite{AlexNick, sequel}. We hope a new candidate  for ``stable $+$ order $+$ random'' will be introduced in the near future.
\end{rem}

\section{Standard models are existentially closed} \label{Sec: Standard models are existentially closed}

\noindent We will finally show that if $\triangleleft$ is a multiplicative circular order on $\ff$, then $(\ff, \triangleleft)$ is a model of $\ACFO$, or in other words, $(\ff, \triangleleft)$ satisfies the geometric characterization. This will require two steps:
\begin{enumerate}
    \item simplifying the characterization of $\ACFO$-models given by Corollary~\ref{refineddefinition of generic} into a characterization that only concerns curves.
    \item using number-theoretic results on character sums over finite fields and counting points over finite fields in combination with Weyl's criterion for equidistribution to show that $(\ff, \triangleleft)$ satisfies the characterization specified in (1).
\end{enumerate}

\subsection{Geometric characterization with curves} \label{sec: Reduction to genericity for curves} We will show that every multiplicatively large variety contains as a subset a curve which is multiplicatively large. This curve will be obtained by intersecting the original variety with suitably chosen hyperplanes of the ambient space. Combining this with Corollary~\ref{refineddefinition of generic}, we will obtain the desired simplified characterization of $\ACFO$-models.

\meno Throughout Section~\ref{sec: Reduction to genericity for curves}, we work with a fixed algebraically closed field $F$, and definable means definable in $F$. The notions of open, closed, irreducible, and dense are with respect to the Zariski topology, which is natural in this context. Let $\dim$ be the canonical dimension for algebraically closed fields, so $\dim$ coincides with Morley rank, topological dimension, \text{acl}-dimension, etc.  Let $\mult$ be the Morley degree. If $X \subseteq F^m$ is definable with $\mult X =1$, we say that $X$ is {\bf generically irreducible} and let the {\bf maximal component} of $X$ be the unique quasi-affine variety with maximal dimension in the decomposition of $X$ into irreducible components.
We let $V$ range over the quasi-affine subvarieties of $(F^\times)^m$ and $C$ range over the one-dimension  quasi-affine subvarieties of $(F^\times)^m$.

\meno Let $S$ be $F^m \setminus \{ 0^{(n)} \}$. If $b$ is an element of $S$, let $H_b$ be the hyperplane defined by  the equation $b\cdot x =1$ where $b \cdot x$ is the usual vector dot product between  $b$ and $x$. So $S$ is essentially the space parametrizing the affine hyperplanes of $F^m$. For each definable set $X \subseteq F^m$  and $b_1, \ldots, b_n \in S$, set
$$  X(b_1, \ldots, b_n)  =  X \cap H_{b_1} \cap \ldots \cap H_{b_n}. $$

\meno Fact~\ref{fact: intersectiondimension} below is a well-known consequence of Fact~\ref{fact: Definable of irreducibility}, Bezout's theorem \cite[Section 4.1]{Shafarevich} and Bertini's theorem \cite[Theorem 2.26]{Shafarevich}.
\begin{fact} \label{fact: intersectiondimension}
Suppose $W \subseteq F^m$ is generically irreducible, and $\dim W = n+1$.  Then the set of $(b_1, \ldots, b_n) \in S^n$ satisfying the following conditions (i) and (ii) is definable and dense in $S^n$:
\begin{enumerate}
    \item[(i)] $W(b_1, \ldots, b_i)$ is generically irreducible for $i \in \{1, \ldots, n\}$;
    \item[(ii)] $\dim  W(b_1, \ldots,b_{i}) = \dim W(b_1, \ldots, b_{i-1}) -1$ for all $i \in \{1,\ldots, n\}.$
\end{enumerate}
Hence, for such $(b_1, \ldots, b_n) \in S^n$, the maximal component of  $W(b_1, \ldots,b_{i})$ is a subset of the maximal component of $W(b_1, \ldots,b_{i-1})$ for $i \in \{1,\ldots, n\}$.
 \end{fact}

\noindent For each $V \subseteq F^m$, define $S_V$ to be the set of $b \in S$ such that $V$ is a subset of $H_b$.

\begin{rem} \label{dimensionofSv}
If $V$ is a single point $c$, then $S_c$ is an irreducible quasi-affine variety and $\dim S_c = m-1$. If $\dim(V) \geq 1$, then $\dim S_V  \leq m-2$.
\end{rem}

\noindent We need a variation of Fact~\ref{fact: intersectiondimension}:

\begin{lem} \label{Lem: variationofBertini}
Suppose $W \subseteq F^m$ is generically irreducible, and $\dim W = n+1$. Then there is $c$ in the maximal component of $W$ such that the set  $Y_c$ of $(b_1, \ldots, b_n) \in S^n_c$ satisfying conditions (i)-(iii) below is dense in $S^n_c$.
\begin{enumerate}
    \item[(i)] $W(b_1, \ldots, b_i)$ is generically irreducible for $i \in \{1, \ldots, n\}$;
    \item[(ii)] $\dim  W(b_1, \ldots,b_{i}) = \dim W(b_1, \ldots, b_{i-1}) -1$ for all $i \in \{1,\ldots, n\};$
    \item[(iii)] $c$ is in $W(b_1, \ldots, b_n)$.
\end{enumerate}
Moreover, with the above $c$, if $X \subseteq W$ is definable and satisfies $\dim X < \dim W$, then the set of  $(b_1, \ldots, b_n) \in Y_c$ such that $\dim X(b_1, \ldots, b_n) \leq 0$   is also dense in $S^n_c$.
\end{lem}
\begin{proof}
 Suppose $W$ and $n$ are as stated above, and $Y \subseteq S^n$ is the set obtained in Fact~\ref{fact: intersectiondimension}.  We show the first statement of the lemma. Let $\Gamma$ be 
$$\{ (b_1, \ldots, b_n,c) \in Y \times W \mid  c \text{ is in the maximal component of  } W(b_1, \ldots, b_n)\}.$$
We note that $\Gamma$ is definable by Fact~\ref{fact: Definable of irreducibility}. Let $\pi_1: \Gamma \to Y$ and $\pi_2: \Gamma \to W$ be the projection maps. For each $c \in W$, we have $\pi_2^{-1}(c) \subseteq S^n_c \times \{c\}$ and $\pi_1( \pi^{-1}_2(c)) \subseteq S^n_c$. We want to find $c$ such that $\pi_1( \pi^{-1}_2(c))$ is dense in $S^n_c$. As $S^n_c$ is irreducible of dimension $n(m-1)$, for the current purpose, it suffices to find  $c \in W$ such that $\dim \pi_2^{-1}(c) = n(m-1)$.

For each $(b_1, \ldots, b_n) \in Y$, we have $ \dim \pi_1^{-1}(b_1, \ldots, b_n) = 1$. As $\dim Y  =  \dim S^n = mn$, it follows that $\dim \Gamma = mn+1$. As $\dim W = n+1$, the set
$$\{ c \in W: \dim \pi_2^{-1}(c) = n(m-1)\}$$
must be dense in $W$. So we obtain $c$ such that the set $Y_c$ of $(b_1, \ldots, b_n) \in S^n_c$ satisfying conditions (i)-(iii) is dense in $S^n_c$.

Let $X$ ne as in the second part of the statement. Note that $\dim X(b_1, \ldots, b_n)$ is at most $\dim W(b_1, \ldots, b_n) =1$ for all $(b_1, \ldots, b_n) \in Y_c$.  By Fact~\ref{fact: Definable of irreducibility} and the irreducibility of $S_c^n$,  exactly one of the following two possibilities happens:
\begin{enumerate}
    \item $\{ (b_1, \ldots, b_n) \in Y_c \mid  \dim X(b_1, \ldots, b_n) \leq 0 \}$ is dense in $S_c^n$;
    \item $\{ (b_1, \ldots, b_n) \in Y_c \mid \dim X(b_1, \ldots, b_n) =1 \}$ is dense in $S_c^n$.
\end{enumerate}
We need to show that (2) cannot happen. Suppose to the contrary that it does. Then using Fact~\ref{fact: Definable of irreducibility}, we get a definable dense subset $R_c$ of $S^n_c$ and  $i \in \{1, \ldots, n\}$ such that 
$$ \dim X(b_1, \ldots, b_{i}) =  \dim X(b_1, \ldots, b_{i-1}) = \dim W(b_1, \ldots, b_{i})  < \dim W(b_1, \ldots, b_{i-1})  $$
for all $(b_1, \ldots, b_n)$  in $R_c$. Shrinking $R_c$ further if necessary, we can arrange that $R_c = U_c^n$ with $U_c$ a definable dense subset of $S_c$. Fix $(b_1, \ldots, b_{i-1}) \in U_c^{i-1}$. Then for all $b_i$ in $U_c$, the hyperplane  $H_{b_i}$ must contain an irreducible component of $ X(b_1, \ldots, b_{i-1}) $ with dimension $\geq 1$. Remark~\ref{dimensionofSv} then gives us that $\dim U_c \leq m-2$ which contradicts the fact that $U_c$ is dense in $S_c$ and $\dim S_c = m-1$. The desired conclusion follows.
\end{proof}

\begin{prop} \label{prop: multiplicative largeness generic intersection}
Suppose $V$ is multiplicatively large and $\dim V = n+1$. Then there is $(b_1, \ldots, b_n) \in S^n$ such that $V(b_1, \ldots, b_n)$ is generically irreducible with  multiplicatively large maximal component of dimension one.
\end{prop}
\begin{proof}
Obtain $c$ in $V$ and $Y_c$ as in the first part of Lemma~\ref{Lem: variationofBertini}. Then for  all $(b_1, \ldots, b_n) \in Y_c$,  $V(b_1, \ldots, b_n)$ is generically irreducible with maximal component $C(b_1, \ldots, b_n)$ of dimension one.
It suffices to show that the set
$$ \{ (b_1, \ldots, b_n) \in Y_c \mid C(b_1, \ldots, b_n) \text{ is multiplicatively large}\} $$
is dense in $Y_c$. Replacing $V$ with $c^{-1}V$ and $c$ with $1^{(m)}$ if neccesary, we arrange that $c = 1^{(m)}$. Hence, $(C(b_1, \ldots, b_m))_{(b_1, \ldots, b_m) \in Y_c}$ is a family of subvarieties of $(F^\times)^m$ passing through $1^{(m)}$. 
Obtain $N$ as in  Corollary~\ref{MultLarge3} for this family. Let $(X_i)_{i=1}^k$ list the intersections of $V$ with the solution sets of the multiplicative equations of the form $t(x)=1$ with  degree at most $N$, and set $X =\bigcup_{i=1}^k X_i.$ 
As $V$ is multiplicatively large, $\dim X < \dim V$.
It then follows from the second part of Lemma~\ref{Lem: variationofBertini} that 
$$ \{ (b_1, \ldots, b_n) \in Y_c \mid  \dim X( b_1, \ldots, b_n) \leq 0 \}  \text{ is dense in } Y_c.$$
  Suppose $(b_1, \ldots, b_n)$ is in the above set. Then  $C(b_1, \ldots b_n)$ is not a subset of $X( b_1, \ldots, b_n)$. So $C(b_1, \ldots b_n)$ is not a subset of $X = \bigcup_{i=1}^k X_i $. By the property of $N$,
$ C(b_1, \ldots, b_n)$ is multiplicatively large, which gives us the desired conclusion.
\end{proof}

\noindent Proposition~\ref{prop: multiplicative largeness generic intersection} together with Corollary~\ref{refineddefinition of generic} yields:

\begin{cor} \label{refineddefinition of generic 2}
Suppose $(F, \triangleleft)$ is a model of $\ACFO^-$ and $(F^\times, \triangleleft) \models \thmm$. Then $(F, \triangleleft)$ satisfies the geometric characterization if and only if all multiplicatively large $C\subseteq (F^\times)^m$ are dense with respect to the $\triangleleft$-topology.
\end{cor}

\subsection{Standard models and number-theoretic randomness} \label{last section}
We now use the results in the preceding section to prove Theorem~\ref{thm: genericity of the standard models}. Other ingredients include a variant of Weyl's criterion for equidistribution and results on counting points and character sums over finite fields, which are consequences of the Weil conjectures for curves over finite fields.

\meno In Section~\ref{last section}, let $\triangleleft$ be the clockwise circular order on $\bbT$.  The multiplicative group $\bbT^m$ is a compact topological group. So $\bbT^m$ is equipped with a unique normalized Haar measure $\mu$. 

\meno A sequence $(X_n)$ of finite subsets of $\bbT^m$ becomes {\bf equidistributed} in $\bbT^m$ if 
$$   \lim_{n \to \infty}  \frac{|X_n \cap U|}{|X_n|} = \mu(U) \quad \text{for all $\triangleleft$-open } U \subseteq \bbT^m . $$
\noindent The following result is a variant of Weyl's criterion for this setting; a proof can be obtained by adapting that of  \cite[Theorem 2.1]{SteinShakarchi}.

\begin{fact} \label{Weylcriterion}
A sequence $(X_n)$ of finite subsets of $\bbT^m$ becomes equidistributed if and only if 
$$   \lim_{n \to \infty} \left( \frac{1}{|X_n|} \sum_{a \in X_n} a_1^{l_1} \cdots a_m^{l_m} \right)=0 \quad \text{for all } (l_1, \ldots, l_m) \in  \zz^m \setminus \{0^{(m)}\}.$$ 
\end{fact}

\noindent Below are the consequences of  Weil conjectures for curves that we need; see \cite{Weilcounting} for Fact~\ref{WeilStyleBound}(i) and \cite[Proposition 4.5]{Perret} for a stronger version of Fact~\ref{WeilStyleBound}(ii).

\begin{fact}\label{WeilStyleBound}
Suppose $C \subseteq \ff^m$ is a one-dimensional quasi-affine variety over $\ff$,  $f \in \ff[C]$ has image in $\ff^\times$, $\text{char}(\ff)=p$,  $C$ and $f$ are definable over $\ff_q$ (in the model-theory sense, or equivalently for perfect fields like $\ff_q$, in the field sense), and $\chi: \ff^\times \to \cc^\times$ is an injective group homomorphism.  Then there is a constant $N \in \nn^{\geq 1}$ such that for all $n\geq 1$, 
\begin{enumerate}
    \item[(i)] $|C(\ff_{q^n})|  < q^n + N \sqrt{q^n}$;
    \item[(ii)] $\big|\sum_{a \in C(\ff_{q^n})} \chi(f(a)) \big|  < N \sqrt{q^n}.$
\end{enumerate}
Here, $C(\ff_{q^n})$ is the set of $\ff_{q^n}$-points of $C$.
\end{fact}

\begin{proof}[Proof of Theorem~\ref{thm: genericity of the standard models}]
Applying Corollary~\ref{Up is model of GMO-3}, it suffices to verify for fixed $\ff$, injective group homomorphism $\chi: \ff^\times \to \bbT$, and multiplicative circular order $\triangleleft_\chi$ on $\ff$ (as defined in the paragraph preceding the same corollary), that $(\ff, \triangleleft_\chi)$ is a model of $\ACFO$.  
By Proposition~\ref{prop: generic and existentially closed}, it suffices to show that $(\ff, \triangleleft_\chi)$ satisfies the geometric characterization. Using Proposition~\ref{Up is model of GMO-} and Corollary~\ref{refineddefinition of generic 2}, we reduce the problem further to showing for a fixed  multiplicatively large one-dimensional quasi-affine variety $C \subseteq (\ff^\times)^m$ that $C$ is  dense in $(\ff^\times)^m$ with respect to the $\triangleleft_\chi$-topology. This is equivalent to showing that $\chi(C)$ is dense in  $\bbT^m$ with respect to the $\triangleleft$-topology.

Assume, without loss of generality, that $\text{char}(\ff)=p$ and $C$ is definable over $\ff_q$. Let $C(\ff_{q^n})$ be the set of $\ff_{q^n}$-points of $C$. Note that $C = \bigcup_n C(\ff_{q^n})$. Hence, the denseness of $\chi(C)$ in $\bbT^m$ with respect to the $\triangleleft$-topology
follows from a stronger result: if $X_n$ is the image of $C(\ff_{q^n})$ under $\chi$, then the  sequence $(X_n)$ becomes equidistributed. Using Fact~\ref{Weylcriterion}, we reduce the problem to verifying that   
$$  \lim_{n \to \infty} \left(  \frac{1}{|C( \ff_{q^n}  )|} \sum_{a \in C( \ff_{q^n}) } \chi (a_1^{l_1} \cdots a_m^{l_m}) \right)=0 \quad \text{for all } (l_1, \ldots, l_m) \in  \zz^m \setminus \{0^{(m)}\}.$$ It is easy to check that $C$ together with $f(x) = x_1^{l_1}\cdots x_m^{l_m}$ satisfies all the conditions described in Fact~\ref{WeilStyleBound},  so we arrive at the desired conclusion.
\end{proof} 

\begin{rem} All approaches to prove Theorem~\ref{thm: genericity of the standard models} so far require the use of character sums over finite fields, counting points over finite fields, and Weyl's criterion for equidistribution. However, slightly different paths could have been taken.

The original approach in our earlier write-up~\cite{Version1} did not go through Section~\ref{sec: Reduction to genericity for curves}, but directly used Corollary~\ref{refineddefinition of generic} and appealed to the much deeper results on character sums over varieties and counting points over varieties \cite{Weil2}. It is possible to avoid these results in the appendix through a rather lengthy proof of ``Lang-Weil Theorem for character sums'' provided in the appendix of the same manuscript.

We also proved there that $\ACFO$ is $\acl_{\text{f}}$-complete (i.e., every complete type over a field-theoretic algebraically closed set is determined by the quantifier-free part of that type); this is a refinement of the fact that $\ACFO$ is model complete as every $\ACFO$-model expands an algebraically closed field. Hrushovski pointed out a shorter path to $\acl_{\text{f}}$-completeness which only uses character sums and counting points over curves: using a similar aproach as in our proof of Theorem~\ref{thm: Model companion of ACFO-}, one can get a similar axiomatization correspoding to the simplified geometric characterization in Corollary~\ref{sec: Reduction to genericity for curves}; then one can show directly that the resulting theory is $\acl_{\text{f}}$-complete by a back-and-forth argument. The reason only results for curves are necessary in this approach roughly comes from the fact that in the back-and-forth argument, one can extend a field-theoretic algebraically closed set each time by an element with transcendence degree $ \leq 1$. 

The current approach is in between. It preserves some of the original intuitive ideas in \cite{Version1} while not appealing to deep number-theoretic results. Proposition~\ref{prop: multiplicative largeness generic intersection} is interesting in its own right, and the technology might be useful elsewhere. The $\acl_{\text{f}}$-completeness of $\ACFO$ can also be obtained now from the general machinery of interpolative fusions \cite{Interpolativefusion}.

\end{rem}

\subsection*{Acknowledgement}
The author would like to thank William Balderrama, Chee-Whye Chin, Lou van den Dries, Ehud Hrushovski, and Erik Walsberg for various discussions and helpful comments on the drafts.
\bibliographystyle{amsplain}
\bibliography{refs}

\end{document}